\documentclass[twoside,reqno]{amsart}
\usepackage{amsmath}
\usepackage{amssymb,latexsym,amstext,amsthm}
\usepackage{verbatim} 
\usepackage{eucal} 
\usepackage{color}
\usepackage[draft]{todonotes}
\usepackage{amsfonts}
\usepackage{mdframed}
\usepackage{lipsum}

\usepackage{pgf,tikz}
\usepackage{mathrsfs}
\usetikzlibrary{arrows}
\definecolor{qqqqff}{rgb}{0.,0.,1.}
\definecolor{xdxdff}{rgb}{0.49019607843137253,0.49019607843137253,1.}

\definecolor{qqqqff}{rgb}{0.,0.,1.}

\newcommand\pref[1]{\textbf{\ref{#1}}}

\newlength\savedwidth 
        \newcommand\whline{\noalign{\global\savedwidth\arrayrulewidth\global\arrayrulewidth 1.2pt}%
        \hline
        \noalign{\global\arrayrulewidth\savedwidth}}

\theoremstyle{plain}
\newtheorem{theorem}{Theorem}[subsection]
\newtheorem{thm}[theorem]{Theorem}
\newtheorem{lem}[theorem]{Lemma}
\newtheorem{cor}[theorem]{Corollary}
\newtheorem{pro}[theorem]{Proposition}

\theoremstyle{definition}
\newtheorem{DEF}[theorem]{Definition}

\newtheorem{exa}[theorem]{Example}

\newtheorem{rem}[theorem]{Remark}

\newtheorem{pgraph}[theorem]{{}}

\newcommand{\sub}{\subseteq}

\newcommand{\la}{Lie algebra }

\newcommand{\ve}{\varepsilon}

\newcommand{\fm}{\frak{m}}

\newcommand{\ep}{\hfill$\Box$}

\def\ad{\hbox{ad}}
\def\andd{\quad\hbox{and}\quad}
\def\sg{\sigma}
\def\a{\alpha}
\def\be{\beta}
\def\te{\theta}
\def\lam{\lambda}
\def\Lam{\Lambda}
\def\ep{\epsilon}
\def\andd{\quad\hbox{and}\quad}

\def\Aut{\hbox{Aut}}

\def\rcross{R^\times}

\def\andd{\quad\hbox{and}\quad}
\def\ind{\hbox{ind}}
\def\vd{\dot{\mathcal V}}

\def\mm{\mathfrak m}
\def\vt{\tilde{\mathcal V}}

\def\a{\alpha}

\def\w{{\mathcal W}}
\def\rd{\dot{R}}

\def\lam{\lambda}
\def\Lam{\Lambda}

\def\ga{\gamma}

\def\la{\langle}
\def\ra{\rangle}

\def\d{\delta}
\def\b{\beta}

\def\sg{\sigma}
\def\bbbz{{\mathbb Z}}

\def\bbbr{{\mathbb R}}
\def\bbbc{\mathbb C}

\def\span{\hbox{span}}

\def\rank{\hbox{rank}}

\def\Aut{\hbox{Aut}}

\def\bbbr{{\mathbb R}}
\def\ep{\epsilon}
\def\fm{(\cdot,\cdot)}
\def\sub{\subseteq}
\def\rds{\dot{R}_{sh}}
\def\rdl{\dot{R}_{lg}}
\def\GL{GL}

\def\u{{\mathcal U}}
\def\V{{\mathcal V}}
\def\w{{\mathcal W}}
\def\b{{\mathcal B}}

\def\mz{\mathcal Z}

\def\rds{\dot{R}_{sh}}
\def\rdl{\dot{R}_{lg}}

\def\rlg{R_{lg}}
\def\rsh{R_{sh}}

\def\mP{\mathcal P}

\def\mar{\mathrm}

\def\mm{\mathcal M}

\def\pp{\mathcal P}
\def\ll{\mathcal L}
\def\hh{\mathcal H}
\def\what{\widehat}
\def\wtilde{\widetilde}

\begin{document}
\setcounter{page}{1}



\author{Saeid Azam$\;^{1}$, Fatemeh Parishani, Shaobin Tan$\;^{2}$}
\title{A characterization of minimal extended  affine root systems\\
(Relations to Elliptic Lie Algebras)}
\address
{Department of Pure Mathematics\\Faculty of Mathematics and Statistics\\
University of Isfahan\\Isfahan, P.O.Box: 81746-73441\\Iran, and\\
School of Mathematics, Institute for
Research in Fundamental Sciences (IPM), P.O. Box: 19395-5746.
} \email{azam@ipm.ir}
\address{Department of Pure Mathematics\\ Faculty of Mathematics and Statistics\\
 University of Isfahan, P.O.Box: 81746-73441, Isfahan, Iran.}
\email{f.parishani93@sci.ui.ac.ir, f.parishani93@yahoo.com}
\address{School of Mathematical Sciences\\Xiamen University\\
 Xiamen, 361005, China.}
 \email{tans@xmu.edu.cn}
 \thanks{$\;^{1}$This work is based upon research funded by Iran National Science Foundation (INSF) under project No.  4001480. The research was in part carried out in the
IPM-Isfahan Branch}
\thanks{$\;^{2}$This research is partially supported by NSF of China (Nos. 12131018, 12161141001)}

\begin{abstract}
Extended affine root systems appear as the root systems of extended affine Lie algebras. A subclass of extended affine root systems, whose elements are called ``minimal" turns out to be of special interest mostly because of the geometric properties of their Weyl groups; they possess the so-called ``presentation by conjugation".
In this work, we characterize minimal extended affine root systems in terms of ``minimal reflectable bases" which resembles the concept of the ``base" for finite and affine root systems. 
As an application, we construct elliptic Lie algebras by means of Serre's type generators and relations.
\end{abstract}
\keywords{Extended affine root systems, extended affine Weyl groups, reflectable bases, elliptic Lie algebras, presentations of Lie algebras}

\subjclass[2010]{17B67 20F55 17B22 17B65}
\maketitle


\setcounter{section}{-1}
\section{Introduction}\label{1}\setcounter{equation}{0}
 \markboth{S. Azam, F. Parishani, S. Tan}{Elliptic Lie algebras}
Motivated by applications in ``the construction of a flat structure for the base space of
the universal deformation of a simple elliptic singularity" K. Saito in \cite{Sa} introduced the concept of an extended affine root system. 
Considering \cite{Sl}, he predicted the ``existence of  Lie algebras
corresponding to {extended affine root systems} which would describe the universal deformation of the simple elliptic singularity". The prediction was partially answered in \cite{H-KT, P}.
Later in 1997, in a systematic study, a class of algebras associated with extended affine root systems was introduced via a set of axioms \cite{AABGP}. Each element of this class is called an ``extended affine Lie algebra". Since then, the theory of extended affine root systems (Lie algebras) (see Definitions \ref{deff1} and \ref{EALA}) and related topics have been under intensive investigation. 
It turns out that despite the striking similarities with the finite, and affine Lie theory, there are enough interesting differences to make the study of extended affine Lie theory a demanding subject. In what follows, we explain the motivation of the present work in two levels; the level of the root system and Weyl group, and the level of Lie algebra.

Level of Root System and Weyl Group: One of the astonishing privileges of finite and affine Weyl groups is that they are Coxeter groups with respect to 
{certain appropriate generating sets.
	In finite types, these generating sets consist of reflections based on elements of a fundamental system of the corresponding root system \cite{C, Hu2, MP}}; thus by \cite[Chapter IV, \S 8, Corollary 3]{B}, these generating sets are ``minimal''.      
In contrast, it is known that the Weyl group of an extended affine root system $R$ is not, in general, a Coxeter group. In fact, one knows that an extended affine Weyl group is a Coxeter group if and only if the corresponding root system is either finite or affine \cite[Theorem 3.6]{Ho}.
Therefore, finding a suitable presentation for an extended affine Weyl group has always been a demanding problem, see for example {\cite{Kr1, Kr2, A1, A2, AS1, AS2, AS3, AS4, AS5, SaT}}. A {possible} approach to achieve such a presentation is as follows. In the theory of root datum, it is known that the Coxeter presentation implies the so-called ``presentation by conjugation'' \cite[Chapter 5,\S 3]{MP}, see Definition \ref{pbc}. Therefore it is natural to ask if an extended affine Weyl group possesses the presentation by conjugation. 
This question was raised first by Y. Krylyuk and was affirmatively answered in \cite{Kr2} for simply laced extended affine Weyl groups of rank$>1$.  Later in a series of works the question was considered for general types, {\cite{A2, Ho, AS3, AS4, AS5}}. In brief, it is revealed that only elements of a subclass of extended affine Weyl groups enjoy the presentation by conjugation. {The root system corresponding to an element of this class is called a ``minimal" extended affine root system} (see Definition \ref{R-min}). For general extended affine Weyl groups, one gives a weaker presentation, called ``the generalized presentation by conjugation'' \cite{A2}, which coincides with the presentation by conjugation if the corresponding root system is minimal. {This explains our motivation for the study of minimal extended affine root systems.}

{The main objective of this work is to give a characterization of minimal extended affine root systems of reduced types in terms of minimal reflectable bases, see Theorem \ref{L11-1}. We achieve this as a by-product of a study of interrelations between three classes of certain subsets of extended affine root systems which we denote by $\mm_r$, $\mm_m$ and $\mm_c$, and we explain them and their related terminologies below.}
Let $R$ be a reduced extended affine root system with the set of non-isotropic roots $R^\times$, and the Weyl group $\w$. For a subset  $\pp$ of  $R^\times$, we denote the subgroup of $\w$ generated by reflections based on $\pp$ by $\w_\pp$. The set $\pp$ is called a {\it reflectable base} for $R$ if $\w_\pp\pp=R^\times$ and no proper subset of $\pp$ has this property. The class of all such sets is denoted by $\mm_r$. The class of all sets $\pp\sub R^\times$ which are minimal with respect to the property that $\w_\pp=\w$ is denoted by $\mm_m$, and  
the class of all sets $\pp\sub R^\times$ which have the minimal cardinality with respect to the property $\w_\pp=\w$ is denoted by $\mm_c$. Let $\mm_{rm}=\mm_r\cap\mm_m$. In the finite theory, one concludes that $\mm_r=\mm_m=\mm_c$, see  Proposition \ref{T3}.

Level of Lie Algebra: K. Saito conjectured in \cite{Sa} that starting from any extended affine root system $R$ there should be a Lie algebra with root system $R$. This led to many interesting investigations
concerning the construction and classification of such Lie algebras, see for example \cite{H-KT, AABGP, Neh}. An extended affine root system of nullity $2$ (see Definition \ref{deff1}) is called an ``elliptic root system'', and the corresponding Lie algebra is called an
``elliptic Lie algebra''. In \cite{SaY, Ya}, the authors construct certain elliptic Lie algebras of rank $>1$ through a Serre's type generators and relations, using a concept of ``base'' which resembles  
the usual notion of base for finite and affine root systems.
The bases used by \cite{SaY, Ya} are not in-general reflectable bases; they involve more roots than reflectable bases. One of our motivations for this work is {to examine the concept of a reflectable base for constructing elliptic Lie algebras by a Serre's type presentation.} We {investigate this} in Section \ref{reala}, in particular, we construct elliptic Lie algebras of rank $1$, the missing case in the works of \cite{SaY, Ya}.   
We now explain the core materials appearing in each section.

In Section 1, we recall from \cite[Definition II.2.1]{AABGP} {the} definition of an extended affine root system and its internal structure in terms of a finite root system and certain subsets of the radical of the form called semilattices, where for our purposes in this work,  we have modified the approach given in \cite[Chapter II, \S 2]{AABGP}, see Proposition \ref{modify1}. 
In Section 2, the Weyl group of an extended affine root system called
an ``extended affine Weyl group'' is considered and a presentation for it called ``the generalized presentation by conjugation'' is explained.  Some immediate related results are derived, see Corollary \ref{cor12} and Lemma \ref{coron-100}. 

{In Section 3, we introduce the classes $\mm_c$, $\mm_r$ and $\mm_m$ for an extended affine root system $R$.
	It is shown that each element $\pp$ of $\mm_r$, $\mm_c$ or $\mm_m$ is a connected generating set for the root lattice, see Proposition \ref{coron-21}. Moreover, associated to each connected subset $\pp\sub R^\times$ a method of assigning an extended affine root subsystem $R_\pp$ is provided such that if $\pp$ belongs to either of $\mm_m$ or $\mm_c$, then $R_\pp$ has the same rank, nullity and type of $R$ and if $\pp\in\mm_r$, $R_\pp=R$, see Proposition \ref{coron2}. Finally, it is derived that for the finite case, 
	we have $\mm_r=\mm_m=\mm_c$ and that this equality fails for general extended affine root systems, see Proposition \ref{T3} and Example \ref{ex1}.}

{In Section 4, {a characterization theorem for reflectable bases associated with reduced extended affine root systems is recorded from \cite[Section 3]{AYY} and
		\cite{ASTY1}. Then} the cardinality of each element $\pp$ in $\mm_r$, $\mm_m$ or
	$\mm_c$ is investigated. It is shown that the cardinality of $\pp$ is finite and, when $R$ is of type $A_1$ or one of the non-simply laced types, all reflectable bases (the elements of $\mm_r$) have the same cardinality that can be described precisely in terms of the rank, nullity, twist number and indices of the involved semilattices. When $R$ is of one of the simply laced types with rank$>1$, $\mm_m=\mm_r$, and reflectable bases may have different cardinalities, see Propositions \ref{L11}, \ref{T2}, \ref{min AD}, and Example \ref{ex2}.}

{In Section 5, the concept of a minimal extended affine root system is recalled from the literature. Starting from a reflectable base $\pp$, a presented group $\widehat\w$ defined by generators ${\widehat w}_{\a}$, $\a\in\pp$, and certain conjugation relations, see Definition \ref{pbc}, is associated to $R$.
	It is shown that $\{\widehat{w}_\a\mid\a\in\pp\}$ is a minimal generating set for $\widehat\w$, see Proposition \ref{min-gen.}. This utilizes proof of the main result of the paper: a reduced extended affine root system $R$ is minimal if and only if $\mm_{r}=\mm_{rm}$; if and only if the corresponding Weyl group has the presentation by conjugation, see Theorem \ref{pro13}. The section is concluded with a table, Table \ref{test}, that illustrates the connections between the results obtained in Sections 1-5.} 

 {We now explain the concluding section, Section \ref{reala}, in which some applications of the concept of a reflectable base are examined at the level of Lie algebra. It is declared that the core $\ll_c$ of an extended affine Lie algebra $\ll$ is finitely generated by showing that it is generated as a Lie algebra by the one-dimensional root spaces $\ll_{\pm\a}$, $\a\in\pp$, $\pp$ being a reflectable base for the ground root system $R$.
	Moreover, if $R$ has index zero (see \pref{note2}(\ref{ind})), then  $\pp$ satisfies a minimality condition concerning this property, Lemma \ref{llem3}.
	The rest of the section deals with constructing an elliptic Lie algebra utilizing a Serre's type presentation based on a reflectable base, see Proposition \ref{proper}.
	As a by-product of this construction, we provide an elliptic Lie algebra of rank $1$ whose non-isotropic root spaces are one-dimensional, the missing case in \cite{SaY, Ya}.} 

{We conclude the introduction by mentioning that in a recent work, see \cite{AFI}, reflectable bases have been used effectively in equipping the core of an extended affine Lie algebra of rank $>1$ with an integral structure, a priority in establishing a modular theory for extended affine Lie algebras.}

\section{Preliminaries}\label{preliminaries}\setcounter{equation}{0}
In this section, we provide some preliminaries on extended affine root systems, in particular, we describe the structure of a reduced extended affine root system in terms of a finite root system and certain subsets called ``semilattices''.  

All vector spaces are finite-dimensional and considered over the field $\bbbr$. For a vector space $\u$, we denote the dual space of $\u$ by $\u^{\ast}$. for a subset $S$ of $\u$, we denote by $\left\langle S\right\rangle$, the additive subgroup of $\u$ generated by $S$. {The symbol $\uplus$ denotes the disjoint union of sets and $|\mP|$  shows the cardinality of a set $\mP$.}

For a vector space $\u$  equipped with a symmetric bilinear form $\fm$, we denote the radical of the form by $\u^0$. For a subset  $R$ of $\u$, we set {$R^0:=R\cap\u^0$} and $R^\times:=R \setminus R^0$. A subset $\pp$ of $\u^\times$ is called {\it connected} (with respect to the form), if
$\pp$ can not be decomposed as $\pp_1\uplus \pp_2$, where $\pp_1$ and $\pp_2$ are non-empty subsets of $\pp$ satisfying $(\pp_1, \pp_2)=\{0\}$.  
\subsection{Extended affine root systems}\label{extended affine root systems}
\begin{DEF}\label{deff1}  Assume that $\V$ is a finite-dimensional real vector space equipped with a positive {semidefinite} symmetric bilinear form $\fm$. Suppose $R$ is a subset of $\V$. Following \cite[Definition 2.1]{AABGP}, we say $R$ is an {\it extended affine root system} if
\begin{enumerate}
\item $0\in R$.
\item $R=-R$.
\item $R$ spans $\V$.
\item If $\a\in R^\times$, then $2\a\not\in R$.
\item $R$ is discrete in $\V$, with respect to the natural topology
of $\V$ when identified with $\bbbr^n$, $n=\dim \V$.
\item If $\a\in R^\times$ and $\be\in R$, then there exist non-negative integers $d, u$ satisfying $\{\be+n\a\mid n\in\bbbz\}\cap R=\{\be-d\a, \ldots, \be+u\a\}$ and $d-u=\frac{2(\a, \be)}{(\a, \a)}$.
\item {$R^\times$ is {connected}}.
\item If $\sg\in R^0$, then there exists $\a\in R^\times$ such that $\a+\sg\in R$.
\end{enumerate}
\end{DEF}

Elements of $R^\times$ are called {\it non-isotropic roots}, and elements of $R^0$ are called {\it isotropic roots}. {The dimension $\nu$ of $\V^0$ is called the {\it nullity} of $R$.}
{Since the form is positive semidefinite, one can easily} check that $$R^0=\{ \alpha \in R\mid (\alpha, \alpha)=0 \}\andd R^{\times}=\{ \alpha \in R\mid (\alpha, \alpha)\neq 0 \}.$$ 

{As it is announced in {\cite{AF}}, the definition of an extended affine root system is equivalent to the following.

\begin{DEF}\label{def2}\rm{
		Let $\V$ be a finite-dimensional real vector space equipped with a positive semidefinite symmetric bilinear form $\fm$, and $R$ a subset of $\V$.  Then $R$ is called an {\it extended affine root system}
		{if the following axioms hold:}}
	{\begin{itemize}
			\item[(R1)] $\la R\ra$ is a full lattice in $\V$.
			\item[(R2)] $(\beta,\a^\vee)\in\bbbz$, $\a,\beta\in R^\times$.
			\item[(R3)] $ w_\a(\beta)\in R$ for $\a\in R^\times$, $\beta\in R$.			
			\item[(R4)] $R^0= \V^0\cap(R^\times - R^\times)$.
			\item [(R5)] $\a\in R^\times\Rightarrow 2\a\not\in R$.
			\item [(R6)] $R^\times$ is connected.
\end{itemize}
}
\end{DEF}
}		

\begin{rem}\label{0rem}
(i) The axiom (R1) means that the natural map
$\la R\ra\otimes_\bbbz \bbbr\rightarrow\V$ is a vector space isomorphism. In other words $\la R\ra$ is the $\bbbz$-span of an $\bbbr$-basis of $\V$. 

(i) From (R4), we conclude that $R^\times\not=\emptyset$.
\end{rem}
{Let $R$ be an extended affine root system in a vector space $\V$ and let $\bar{\;}:\V\rightarrow\bar\V:=\V/\V^0$ be the canonical map.
 It is known that $\bar{R}$, the image of $R$ under the map $\bar{\;}$ is {an irreducible finite root system}  in $\bar\V$, see \cite[Proposition II.2.9(d)]{AABGP}. 
 
 \begin{DEF}\label{0DEF} 
The  {\it type} and the {\it rank} of extended affine root system $R$ is defined as the type and the rank of $\bar R$. In this work, we always assume that $R$ is of {\it reduced type}, that is $\bar R$ has one of the types $A$, $B$, $C$, $D$, $E_{6,7,8}$, $F_4$ or $G_2$. 
\end{DEF}

 %

\begin{DEF}\label{1def}
(i) For a real vector space $\u$ equipped with a symmetric bilinear form $\fm$, and $\a\in\u$ with $(\a,\a)\neq 0$, the {\it reflection based on} $\a$ is defined by  
 $\be\mapsto\be-\frac{2(\be,\a)}{(\a,\a)}{\a},\;\be\in\u.$
 By convention, if $(\a,\a)=0$, we interpret $w_\a$ as the identity map.
 
 (ii) For a subset $\pp\sub\u$, we denote by $\w_\pp$ the subgroup of $\GL(\u)$ generated by $w_\a$, $\a\in\pp$.
 \end{DEF}

\subsection{Semilattices}
Let $R$ be an extended affine root system in $\V$ with corresponding from
$\fm$. As was mentioned before, we denote the radical of $\fm$ by $\V^0$. Certain subsets of $\V^0$ called ``semilattices'' play an important role in describing the internal structure
of extended affine root systems. 
\begin{DEF}\label{0sem}
	A {\it semilattice} in $\V^0$ is a spanning subset $S$ of $\V^0$ which is discrete, contains $0$ and satisfies $S\pm 2S\sub S$.
\end{DEF}

Here we record some facts about semilattices needed in the sequel.
\begin{pgraph}\label{not1}
	Let $S$ be a semilattice  
	in $\V^0$ and set $\Lam:=\la S\ra$.
	
	(i)  $\Lam$ is a lattice in $\V^0$; that is $\Lam$ is the $\bbbz$-span of a basis of $\V^0$.
	
	(ii) $2\Lambda+ S\subseteq S\sub\Lam$ and so
	$S=\bigcup_{i=0}^{m}(\tau_i +2\Lam),$ for some
	$\tau_i\in S$, with	$\tau_0=0$.

	(iii) $\Lam$ admits a $\bbbz$-basis $B=\{\sg_1,\ldots,\sg_\nu\}$ consisting of elements of $S$. 
\end{pgraph}

\begin{DEF}\label{-1def} Assume that {$S$ is a semilattice}
	with $\Lam=\la S\ra$.
	
	(i) The {\it index} of  $S$, denoted $\ind(S)$, is by definition the least positive integer $m$ such that
	\pref{not1}(ii) holds.
	
	(ii) Assume that $\ind(S)=m$, and $\tau_i$'s and  $B$ are as in \pref{not1}(iii).
	For $\sg=\sum_{i=1}^\nu n_i\sg_i\in\Lam$, we define $\mbox{supp}(\sg):=\{i\mid n_i\in 2\bbbz+1\}.$
	Then $\sg=\sum_{i\in\mbox{supp}(\sg)}\sg_i$ (mod $2\Lam$).
	The collection
	$\mbox{supp}(S)=\{\mbox{supp}(\tau_i) \mid  0\leq i\leq m\},$
	is called the {\it supporting class} of $S$ with respect to the basis $B$.
\end{DEF}

The supporting class determines $S$ uniquely. In fact, we have
\begin{equation}
S=\biguplus_{J\in \mbox{supp(S)}}(\tau_J+2\Lam),\label{e2}
\end{equation}
where $\tau_J:=\sum_{r\in J}\sg_r$.

\subsection{Internal structure of extended affine root systems}\label{Internal structure}
In \cite[Chapter II.\S2]{AABGP}, the integral structure of an extended affine root system in terms of a finite root system and certain semilattices is studied. Based on our purposes in this work, we need to make some modifications to \cite{AABGP}'s approach.
As before we assume that $R$ is a reduced extended affine root system in $\V$ equipped with the positive semi-definite form $\fm$. 

We start with recording the following lemma whose proof for simply laced types is given in \cite[Proposition 2.8]{ASTY1}, and for non-simply laced types concludes from \cite[Propositions 2.6, 2.9 and Remark 2.15]{AYY}. 

\begin{lem}\label{0lem}
	Let $\Phi$ be a reduced irreducible finite root system in a vector space $\u$, and $\pp$ be a minimal subset (with respect to inclusion) of $\Phi^\times$ such that
	$\w_\pp\pp=\Phi^\times$. Then $|\pp|=\rank\; \Phi$.
\end{lem}


\begin{pro}\label{modify1}
	Let $R$ be an extended affine root system and $\pp$ a subset of  $R^\times$ satisfying $\w_{\pp}\pp=R^\times$. Then
	$\pp$ contains a subset $\dot\Pi$ such that the following hold.
	
	(i)  $\dot R:=(\w_{\dot\Pi}\dot\Pi)\cup\{0\}$ is a finite root system
	in $\dot\V:=\span_{\bbbr}\dot R$,
	
	(ii) $\dot R$ is isomorphic to $\bar R$,
	
	(iii) $\V=\dot\V\oplus\V^0,$
	
	(iv) For a fixed $\a\in R^\times$, there exists a finite root system ${\dot R}_\a\sub R$  such that $\a\in{\dot R}_\a$ and (ii) and (iii) hold with  ${\dot R}_\a$ and ${\dot\V}_\a:=\span_{\bbbr}\rd_\a$ in place of $\rd$ and $\dot\V$.  
 \end{pro}
 
\proof {Since $R^\times=\w_\pp\pp$, we have $\bar{R}^\times=\w_{\bar\pp}\bar\pp$. Let $\bar\Pi$ be a subset of $\bar\pp$ such that $\bar{R}^\times=\w_{\bar\Pi}\bar\Pi$ and no proper subset of $\bar\Pi$ satisfies this property. From Lemma \ref{0lem}, we have 
 $|\bar\Pi|=\rank\;{\bar R}$ and so $\bar\Pi$ is a basis of the vector space $\bar\V$.  Let $\bar\Pi=\{\bar{\a}_1,\ldots,\bar{\a}_\ell\}$. We fix a preimage ${\dot\a}_i \in\pp$ for  ${\bar\a}_i$. Set $\dot\Pi:=\{{\dot\a}_1,\ldots,{\dot\a}_\ell\}$, $\dot\V:=\span_{\bbbr}\dot\Pi$, and
 $\dot R:=(\w_{\dot\Pi}\dot\Pi)\cup\{0\}$. Then the epimorphism $\bar{\;}:\V\rightarrow\bar\V$ induces an isometry $\dot\V\rightarrow\bar\V$ which maps $\dot R$ isometrically onto $\bar R$. It follows that $\V=\dot\V\oplus\bar\V$. It completes the proof of (i)-(iii).}

 {Next, let $\a\in R^\times.$ Then $\bar\a\in\bar R=\w_{\bar\Pi}\bar\Pi$, and so $\bar\a=w_{{\bar\a}_{i_1}}\cdots w_{{\bar\a}_{i_k}}({\bar\a}_{i_{k+1}})$ for some ${\bar\a}_{i_1},\ldots, {\bar{\a}}_{i_k}, {\bar{\a}}_{i_{k+1}}\in\bar\Pi$. Set 
 $$
  {\bar{\Pi}_\a}:=w_{{\bar\a}_{i_1}}\cdots w_{{\bar\a}_{i_k}}({\bar\Pi}).
  $$
Now ${\bar\Pi}_\a$ is a basis of $\bar\V$ such that $\bar\a\in{\bar\Pi}_\a$ and $\w_{{\bar\Pi}_\a}{\bar\Pi}_\a={\bar R}^\times$. Therefore,  in the proof of part (i), we may work with ${\bar{\Pi}}_\a$ in place of ${\bar\Pi}$. Then as above we fix a preimage
${\dot\Pi}_\a$ of {${\bar\Pi}_\a$} where we pick $\a$ as the preimage of $\bar\a$. Finally, we construct ${\dot R}_\a$ in the same manner as $\dot R$. Then ${\dot R}_\a$ and ${\dot\V}_\a:=\span_{\bbbr}{\dot R}_\a$ satisfy (ii)-(iii) in place of
$\dot R$ and $\dot\V$, and the proof of (iv) is completed.
 \qed
 }

 Let $\rd$ and $\vd$ be as in Proposition \ref{modify1}. 
We write $\rd^\times=\rds\cup\rdl$ where $\rds$ and $\rdl$ are the sets of short and long roots of $\rd$. 
 By convention, we assume that all non-zero roots of $\rd$ are short if there is only one root length in $\rd$.

\begin{pgraph}\label{pgraph2}
It is shown in \cite[Chapter II]{AABGP} that using the finite root system $\dot R$, one can obtains a description of $R$ in the form
\begin{equation}\label{e11}
\begin{array}{l}
R=(S+S)\cup(\rds+S)\cup (\rdl+L),\hbox{ where}\\
S=\{\sg\in\V^0\mid\sg+\a\in R\hbox{ for some }\a\in\rds\},\\
L=\{\sg\in\V^0\mid\sg+\a\in R\hbox{ for some }\a\in\rdl\}.
\end{array}
\end{equation}
where $S$ and $L$ are {\it semilattices} in $\V^0$.
If $\rdl=\emptyset$,  we interpret the terms $L$ and  $ \rdl+L$ as  empty sets. Semilattices $S$ and $L$ satisfy some further properties which we record here. 
\end{pgraph}

\begin{pgraph}\label{not2}
If $R$ is simply laced of rank $>1$, or is of type $C_\ell$, $\ell\geq 3$, $F_4$ or $G_2$, then $S$ is a lattice. If $R$ is of type $B_\ell$, $\ell\geq 3$, $F_4$ or $G_2$, then $L$ is a lattice. Furthermore, for non-simply laced types, we have
\begin{equation}\label{e1-1}
L+S=S\andd k S +L=L ,
\end{equation}
with 
\begin{equation}\label{k23}
k=\left\lbrace\begin{array}{ll}
3&X=G_2,\\
2&\hbox{other types.}
\end{array}\right.
\end{equation}
 Now (\ref{e1-1}) implies that 
\begin{equation}\label{e12}
k\la S\ra\subseteq \la L\ra\subseteq S;
\end{equation} 
and so $\la S\rangle/{\la L\ra}$ is a finite-dimensional vector space over the finite field $\bbbz_k$. The integer $0\leq t\leq\nu$ satisfying $|{\la S\rangle}/{\la L\ra}|=k^t$ is called the {\it twist number} of $R$. {We make the convention that if $R$ is of simply laced type, then $t=0$ and $k=1$.  } 
\end{pgraph}

We proceed with a fixed description of $R$ in the form (\ref{e11}). We assume that $R$ is of type $X$, has rank $\ell$, nullity $\nu$ and twist number $t$. 

\begin{pgraph}\label{note2}
From \cite[II. \S 4(b)]{AABGP}, we see that  there exist subspaces $\V_1^0$ and $\V_2^0$ of vector space $\V^0$ of dimensions $t$ and $\nu -t$ respectively such that
$\V^0 =\V_1^0\oplus\V_2^0,$
and there are semilattices $S_1$ and $S_2$ in $\V_1^0$ and $\V_2^0$, respectively such that 
\begin{equation}\label{s1}
S= S_1\oplus\la S_2\ra\andd
L=k\la S _ 1\ra\oplus S _ 2.
\end{equation}
\end{pgraph}
\noindent
Following \cite[Section 4]{A1}, we define the {\it index} of $R$, denoted $\ind(R)$, by 
\begin{equation}\label{ind}
\mbox{ind}(R):=\left\lbrace \begin{array}{ll}
0& X= A_\ell(\ell\geq 2), D_\ell(\ell\geq 4), E_{6, 7, 8}, F_4, G_2,\\
\mbox{ind}(S) -\nu & X= A_1,\\
\mbox{ind}(S_1)+\mbox{ind}(S_2)-\nu & X= B_2,\\
\mbox{ind}(S_1)- t & X= B_\ell(\ell\geq 3),\\
\mbox{ind}(S_2)- (\nu - t) & X= C_\ell(\ell\geq 3).
\end{array}
\right.
\end{equation}
We write $S_1$ and $S_2$ in the form \pref{not1}(ii), namely 
\begin{equation}\label{coron6}
S_1=\biguplus_{i=0}^{\ind (S_1)}(\gamma_i+2\la S_1\ra)\andd S_2=\biguplus_{i=0}^{\ind (S_2)}(\eta_i+2\la S_2\ra).
\end{equation}

\begin{pgraph}\label{pgraph1}
We fix a fundamental system $\dot{\Pi}=\{\a_1, \ldots, \a_\ell\}$ for $\rd$. By \pref{not1}(iii), $\Lam=\la S\ra$ admits a $\bbbz$-basis $\{\sg_1,\ldots,\sg_\nu\}\sub S$. Let $\theta$ be any root in $\rd$ {if} $\rd$ is of simply laced type, and $\te_s$ and $\te_l$ be any short and any long root in $\dot{R}$, respectively if $\rd$ is of non-simply laced type. For later use, for each type $X$, we introduce a subset $\pp(X)$ of $R^\times$ as follows:
\begin{table}[htb]
	\caption{The sets $\pp(X)$}
	\label{tab:fq}
	\vspace{-4mm}
	\scriptsize{
		\begin{tabular}
			[c]{|c |   c |} \hline
					$\mbox{Type}$& $\pp(X)$  \\
					\whline
					$A _1$ & ${\{\a_1, \tau_1 -\a_1, \ldots, \tau_{\rho}-\a_1 \}},\;(\rho=\ind(S))$  \\
					\hline
					$A_\ell (\ell>1)$ & $\{\a_1, \ldots, \a_\ell, \sg_1-\te, \ldots, \sg_\nu-\te\}$\\
					\hline 
					$D_\ell$ & $\{\a_1, \ldots, \a_\ell, \sg_1-\te, \ldots, \sg_\nu-\te\}$ \\
					\hline
					$E_{6,7,8}$ & $\{\a_1, \ldots, \a_\ell, \sg_1-\te, \ldots, \sg_\nu-\te\}\; {(\ell=6, 7, 8)}$\\
					\hline
					$F_4$ & $\{{\a_1, \a_2, \a_3, \a_4}, \sg_1-\te_s, \ldots, \sg_t-\te_s, \sg_{t+1}-\te_l,
					\ldots, \sg_\nu -\te_l\}$ \\
					\hline
					$G_2$ & $\{{\a_1, \a_2}, \sg_1-\te_s, \ldots, \sg_t-\te_s, \sg_{t+1}-\te_l,
					\ldots, \sg_\nu -\te_l\}$\\
					\hline
					$B_2$ & $\{\a_1, \a_2, \gamma_1-\te_s, \ldots, \gamma_{\rho_1}-\te_s, 
					\eta_1 -\te_l, \ldots, \eta_{\rho_2}-\te_l\},\;(\rho_i=\ind(S_i),\;i=1,2)$\\
					\hline
					$B_\ell\; (\ell>2)$ & $\{\a_1, \ldots, \a_\ell, \gamma_1-\te_s, \ldots, \gamma_{\rho_1}-\te_s, \sg_{t+1}-\te_l, \ldots, \sg_\nu -\te_l \},\;(\rho_1=\ind(S_1))$\\
					\hline
					$C_\ell \;(\ell>2)$ & $\{\a_1, \ldots, \a_\ell, \sg_1 -\te_s, \ldots, \sg_t -\te_s, \eta_1 -\te_l, \ldots, \eta_{\rho_1}-\te_l\},\;(\rho_2=\ind(S_2))$\\
					\hline
				\end{tabular}
			}
		\end{table}
	\end{pgraph}

\section{Extended affine Weyl groups }\label{extended affine Weyl groups}\setcounter{equation}{0}
We proceed with the same notation as 
in Section \ref{preliminaries}, in particular we assume that $R$ is a reduced extended affine root system in $\V$ of type $X$, rank $\ell$ and nullity $\nu$. As before, $\V^0$  is the radical of the form.
\subsection{The setting}\label{the setting}
\begin{pgraph}\label{e1}
	We fix a description of $R$ in the form $R=(S+S)\cup(\rds+S)\cup(\rdl+L)$ as in \pref{pgraph2}. We recall that $\dot\V:=\span_\bbbr\rd$ and $\V=\dot\V\oplus\V^0$.  We consider the $\ell+2\nu$-dimensional real vector space $\tilde{\V}:= \dot{\V}\oplus \V^0 \oplus (\V^0)^\ast$.
	We normalize the form on $\dot{\V}$ such that
	$$
	(\dot\alpha,\dot \alpha):=2\;\;\hbox{for}\;\;\dot{\a}\in\rds.$$
	We extend the form $\fm$ on $\V$ to a symmetric form on $\tilde{\V}$, denoted again by $\fm$,  as follows:
	$$\begin{array}{c}
		(\dot{\V}, (\V^0)^\ast)= ((\V^0)^\ast, (\V^0)^\ast):=\{0\}, \\
		(\sg,\lam)=\lam(\sg),\hbox{ for }\sg\in \V^0,\;\lam\in (V^0)^\star.
	\end{array}$$
	This forces the form on $\tilde{\V}$ to be non-degenerate. For $\a\in\V$ with $(\a,\a)\neq 0$ set $\a^\vee:=2\a/(\a,\a)$. 
	
	\begin{DEF}\label{eawg}
		The {\it extended affine Weyl group} 
		$\w$ of $R$ (or the Weyl group of extended affine root system $R$) is defined to be the subgroup of $\GL(\vt)$ generated by reflections $w_\a:\beta\mapsto\beta-(\beta,\a^\vee)\a$, $\a\in R^\times$. 
		Since $0\in L\sub S$, we have $\rd\sub R$, and so we may identify $\dot\w$, the Weyl group of $\rd$, as a subgroup of $\w$. 
		We denote by $\mz(\w)$, the center of $\w$.
	\end{DEF}
\end{pgraph}

\subsection{Presentation by conjugation}\label{presentation by conjugation}
It is known that an extended affine Weyl group $\w$ associated with $R$ is not, in general, a Coxeter group. It is proved that $\w$ is a Coxeter group if and only if $R$ has nullity $\leq 1$ (see \cite[Theorem 3.6]{Ho}).
{Nevertheless, extended affine Weyl groups enjoy another interesting presentation called ``generalized presentation by conjugation''. 
	When $R$ has some ``minimality condition'', then $\w$ admits even a more delightful presentation, called the ``presentation by conjugation''. We record these definitions here, as they play a crucial role in the sequel.}

\begin{pgraph}\label{pgraph3}
	{We begin by recalling some notation from \cite[Sections  1, 2]{A2}. We fix the $\bbbz$-basis $\{\sg_1, \ldots, \sg_\nu\}\sub S$ of $\Lam=\la S\ra$ and the integer $k$ given in \pref{not1}(iii) and \pref{not2}. For $1\leq i\leq\nu$ and $\a\in\la R\ra$ define the linear map $t_{\a}^{(i)}:\tilde{\V}\rightarrow\tilde{\V}$ by 
		$$t_{\a}^{(i)}(\lam):=\lam -(\lam, \frac{1}{k}\sg_i)\a +((\lam , \a)-\frac{1}{2}(\a, \a)(\lam, \frac{1}{k})\sg_i)\frac{1}{k}\sg_i.$$}
	The map $t^{(i)}_\a$ is invertible with the inverse $t^{(i)}_{-\a}$, so
	$t^{(i)}_\a\in\Aut(\tilde{\V}).$  {For $1\leq i, j\leq\nu$ set 
		$$c_{ij}:=t_{-\sg_j}^{(i)}.$$
		One  checks that $(t_{\a}^{(i)}, t_{\be}^{(j)})=c_{ij}^{\frac{1}{k}(\a, \be)}$, for $\a, \be\in R$, where $(x, y):=xyx^{-1}y^{-1}$ denotes the group commutator.}
\end{pgraph}

{ 
For $\a\in\rcross$ and $\sg=\sum_{i=1}^{\nu}m_i\sg_i\in R^0$ with $\a +\sg\in R$, we set
\begin{equation}\label{gpbc}
c_{(\a,\sg)}:=(w_{\a +\sg}w_\a)(w_\a w_{\a +\sg_1})^{m_1}\ldots (w_\a w_{\a +\sg_\nu})^{m_\nu}.
\end{equation}
By \cite[(2.4)]{A2}, we have
\begin{equation}\label{cij}
c_{(\a, \sg)}=\prod_{i<j}c_{ij}^{k(\a)m_im_j},
\end{equation}
where the integers $k(\a)$ are defined as follows.
$$k(\a):=\left\lbrace
\begin{array}{lll}
	1&\hbox{if}&\a\in\rlg,\\
	k&\hbox{if}&\a\in\rsh,
	\end{array}\right.$$
(If $R$ is simply laced, then $k(\a)=1$ for all $\a\in\rcross$).}

\begin{pgraph}\label{pgraph4}
{Assume that $\theta_s$ and $\theta_l$ are the highest short and highest long roots of $\rd$ respectively (for simply laced cases $\theta:=\theta_s=\theta_\ell$), and for $1\leq p\leq\nu$  consider a triple 
$$(\ve_p, \a_p, \eta_p=\sum_{i=1}^{m}m_{ip}\sg_i)\in (\{\pm 1\}\times\{\theta_s\}\times\sum_{i=1}^{\nu}\bbbz\sg_i)\cup(\{\pm 1\}\times\{\theta_l\}\times\sum_{i=t+1}^{\nu}\bbbz\sg_i).$$}
{The collection $\{(\ve_p, \a_p, \eta_p)\}_{p=1}^{m}$ is called a {\it reduced collection} if 
$$\sum_{p=1}^{m}k(\a_p)\ve_pm_{ip}m_{jp}=0~~~\hbox{for all}~~1\leq i<j\leq\nu.$$
 If $\{(\ve_p, \a_p, \eta_p)\}_{p=1}^{m}$ is a reduced collection, then $c_{(\a_p, \eta_p)}$ is an element of $\w$ for all $1\leq p\leq m$ (see \cite[Section 1]{A2}). } 
\end{pgraph}

{The following theorem is due to \cite[Theorem 2.7]{A2}, \cite[Theorem 1.18]{Kr1}, {\cite[Theorem 5.19(ii)]{AS3}} and \cite[Corollaries 2.34, 2.36(b)]{AS5}.}
 \begin{thm}\label{pbc1}
Let $R$ be an extended affine root system  with Weyl group $\w$ and $\widehat{\w}$ be the group  {defined by generators} $\widehat{w}_\a$, $\a\in R^\times$, subject to {the} relations:
\begin{equation}\label{(i)}\widehat{w}_\a^2=1\hbox{ for every }\a\in R^\times,\vspace{-3mm}\end{equation}
\begin{equation}\label{(ii)}
\widehat{w}_\a\widehat{w}_\be\widehat{w}_\a=\widehat{w}_{w_\a(\be)}\hbox{ for every }\a, \be\in R^\times,\vspace{-3mm}
\end{equation} 
\begin{equation}\label{(iii)}
\prod_{p=1}^{m}\widehat{c}^{\ve_p}_{(\a_p, \eta_p)}=1\hbox{ for any reduced collection }\{(\ve_p, \a_p, \eta_p)\}_{p=1}^{m},\end{equation}
where $\widehat{c}_{(\a_p, \eta_p)}$ is the element in $\widehat{\w}$ corresponding to $c_{(\a_p, \eta_p)},$ under the assignment $w_\a\mapsto \widehat{w}_\a.$

(i) The extended affine Weyl group $\w$ is isomorphic to  $\widehat{\w}.$ 

(ii) If $R$ is simply laced of rank$>1$, or is of type $F_4$ or $G_2$ then relations of the form (\ref{(iii)}) are consequences of relations of the forms (\ref{(i)}) and (\ref{(ii)}); in the sense of Definition \ref{pbc} below, $\w$ has the presentation by conjugation.

(iii) If $\nu\leq 2$, then $\w$ has the presentation by conjugation (Definition \ref{pbc}).

\end{thm}

\begin{DEF}\label{pbc}
	The given presentation in Theorem \ref{pbc1} is called the {\it generalized presentation by conjugation}. 
	If $\w$ is isomorphic to the presented group defined by generators
	$\widehat{w}_\a$, $\a\in R^\times$ and relations (\ref{(i)}) and (\ref{(ii)}), then $\w$ is said to have the {\it presentation by conjugation.} \end{DEF}

\begin{pgraph}\label{pgraph7}
	As a consequence of \cite[Theorems 4.23, 5.16]{AS3} and \cite[Theorem 2.33]{AS5},  the extended affine Weyl group $\w$ has the presentation by conjugation if and only if the epimorphism
	\begin{eqnarray}\label{epi}
		\psi:\widehat{\w}\longrightarrow\w\\
		{\widehat{w}}_\a\mapsto w_\a,\nonumber
	\end{eqnarray}
	is an isomorphism.
\end{pgraph}

{For $\a\in R^\times$, we denote the orbit of $\a$ under the action of $\w$ on $R^\times$ by $\w\a$, namely $\w\a=\{w\a\mid w\in\w\}$.}
\begin{cor}\label{cor12}
	Assume that $\w$ has the presentation by conjugation and fix $\beta\in R^\times$. Then the assignments 
	$$
	w_\a\stackrel{\Phi_\be}{\longmapsto}\left\{\begin{array}{ll}
		1 & \mbox{if}~ \a\in\w\be,\\
		0 & \mbox{if} ~\a\in\rcross\setminus\w\be,
	\end{array}\right.\quad\quad\quad
	w_\a\stackrel{\Psi}{\longmapsto}\left\{\begin{array}{ll}
		0&\hbox{if }\a\in R_{lg},\\
		1&\hbox{if }\a\in R_{sh}.
	\end{array}\right.
	$$
	extend to epimorphisms
	$ \Phi_\be$ and  ${\Psi}$ from $\w$ onto $\bbbz_2$,  respectively.
\end{cor}

\proof First, we consider $\Phi_\beta$. Since for $\a,\gamma\in R^\times$, we have
$\gamma\in\w\beta$ if and only if $w_\a(\gamma)\in\w\beta$, it easily follows that the defining relations of the form 
$w_\a w_\gamma w_\a w_{w_\a(\gamma)}$ are mapped to zero under the assignment $\Phi_\beta$. The argument for defining relations of the form  $w_\a w_\a$ is obvious. 

For the assignment $\Psi$, clearly, the defining relations are mapped to zero under $\Psi$,  as reflections preserve the length.\qed

	
	\begin{lem}\label{coron-100}
		{Let $R$ be of non-simply laced type. Let $\{w_\a\mid\a\in\pp\}$
			be a set of generators for $\w$. Then $\pp\cap R_{sh}\neq\emptyset$ and $\pp\cap R_{lg}\neq\emptyset$.}
	\end{lem}
	
	\proof 
	First, assume that $\w$ has the presentation by conjugation and consider the epimorphism $\Psi:\w\rightarrow\bbbz_2$ of Corollary \ref{cor12}.  Now if $\pp\sub R_{\lg}$, then for $\a\in R_{sh}$, we have
	$w_\a=w_{\a_1}\cdots w_{\a_m}$ for some $\a_1,\ldots,\a_m\in R_{lg}$. Then
	$1=\Psi(w_\a)=\Psi(w_{\a_1})+\cdots+\Psi(w_{\a_m})=0$ which is absurd. An analogous argument works for the case $\pp\sub R_{sh}$, replacing the roles of short and long roots. Thus the result holds if $\w$ has the presentation by conjugation.} 

{For the general case, we consider the epimorphism
	$\w\rightarrow\w_{\bar R}$, $w_\a\mapsto w_{\bar\a}$, where
	$\w_{\bar R}$ denotes the Weyl group of the finite root system $\bar R$. Now since $\{w_\a\mid\a\in\pp\}$ generates $\w$, the set
	$\{w_{\bar\a}\mid\a\in\pp\}$ generates $\w_{\bar R}$. Since $\w_{\bar R}$ has the presentation by conjugation {(see \cite[Proposition 5.3.3]{MP})}, we get from the first part of the proof that ${\bar R}_{sh}\cap\bar\pp\not=\emptyset$ and 
	${\bar R}_{lg}\cap\bar\pp\not=\emptyset$.
	Since for $\a\in R$, $\a$, and $\bar\a$ have the same length, we are done.}\qed

We return to the discussion about the existence of presentation by conjugation for extended affine Weyl groups in {future sections.}

\section{The classes $\mm_r$, $\mm_c$ and $\mm_m$}\setcounter{equation}{0}
We proceed by assuming that  $R$ is a reduced extended affine root system in $\V$ of rank $\ell$ and  nullity $\nu\geq 0$. Let $\w$ be the Weyl group of $R$.
We begin with introducing some notation.
For a subset $\pp$ of $R^\times$, we set
$$\begin{array}{l}
\pp_{sh}:=\pp\cap\rsh,\\
\pp_{lg}:=\pp\cap\rlg,\\
S_\pp=\{w_\a\mid\a\in\pp\},\\
\w_\pp:=\la S_\pp\ra=\la w_\a\mid\a\in\pp\ra.
\end{array}$$

\subsection{Inter relations between  $\mm_r$, $\mm_c$ and $\mm_m$}\label{reflectable bases}
{As we will see} in the sequel (Proposition \ref{L11}), the extended affine Weyl group $\w$ is finitely generated. 
\begin{DEF} \label{D1}
Let $c$ denote the least number of reflections generating $\w$.
Let $\mP\subseteq R^\times$; we say $\mP$ is a {\it $c$-minimal} set {in $R^\times$} if, $|\pp|=c$ and $S_\pp$ generates $\w$. 
%
%
\end{DEF}

\begin{DEF}\label{D11}
(\cite[Definition 1.19]{AYY})  
A subset $\pp\subseteq R^\times$ is called a {\it reflectable set}
if $\w_\pp\pp=R^\times$. A subset $\pp\sub R^\times$  is called a {\it reflectable base} if

(i) $\pp$ is a reflectable set; $\w_\pp\pp=R^\times$,

(ii) no proper subset of $\pp$ is a reflectable set. 
\end{DEF}
\begin{DEF}\label{D2} 
	Let $\mP\subseteq R^\times$; $\mP$ is called a {S-\it minimal} set (or simply a {\it minimal} set) {in $R^\times$} if,
	
	{(i)} $S_\mP$ generates $\w$,
	
	{(ii)} no proper subset of $S_\mP$ generates $\w$. 
	
\end{DEF} 

\begin{pgraph}\label{pgraph9}
{We denote the class of reflectable bases, $c$-minimal sets, and $S$-minimal
	sets of $R^\times$ by ${\mm}_r$, $\mm_{c}$ and $\mm_m$, respectively. We also set 
	$$\mm_{rc}:=\mm_r\cap\mm_{c}\andd
	\mm_{rm}:=\mm_r\cap\mm_{m}.
	$$}
Clearly, we have
\begin{equation}\label{coron1}
	\mm_c\sub\mm_m\quad\hbox{and }\quad\mm_{rc}\sub\mm_{rm}.
\end{equation}
{Note that if $\pp\in\mm_r\cup\mm_m$, then as $\w=\w_\pp$, we have $\V=\span_{\bbbr}\pp$, and then 
	\begin{equation}\label{coron11}
		\dim(\V)=\ell +\nu\leq c\leq |\pp|.
\end{equation}}
\end{pgraph}

\begin{pro}\label{coron-21}
If $\pp\sub R^\times$ satisfies $\w=\w_\pp$, then $\la \pp\ra=\la R\ra$ and $\pp$ is connected. {In particular, if $\pp$ belongs to either of  $\mm_c$, $\mm_m$ and $\mm_r$, then $\la\pp\ra=\la R\ra$ and $\pp$ is connected.}
\end{pro}

\proof {For $\a\in R^\times$, we have $w_\a=w_{\a_1}\cdots w_{\a_n}$ for some $\a_1,\ldots,\a_n\in\pp$. If there exists $\beta\in R^\times$ with $(\beta,\a^\vee)=-1$, then
$$\a+\beta=w_\a(\beta)=w_{\a_1}\cdots w_{\a_n}(\beta)=\beta+k_1\a_1+\cdots+k_n\a_n$$
for some $k_1,\ldots,k_n\in\bbbz$.  Thus $\a\in\la\pp\ra$ and so $\la\pp\ra=\la R\ra$. In particular, this holds for  the types
$X\not=A_1,\; B_\ell,\;C_\ell$. Assume now that $R$ is of type $A_1$. Then
$$-\a=w_\a(\a)=w_{\a_1}\cdots w_{\a_n}(\a)=\a+2k_1\a_1+\cdots+2k_1\a_n,$$ for some $k_i\in\bbbz$, which again give $\a\in\la\pp\ra$, as required. }

{Next assume that $R$ is of type $B_\ell$.  If $\a\in R_{lg}$, then there exists $\beta\in R_{sh}$,
such that $(\beta,\a^\vee)=-1$, and so we have $\a\in\la\pp\ra$ by the preceding paragraph. 
Thus $\la R_{lg}\ra\sub\la\pp\ra.$ It remains to show that $R_{sh}\sub\la\pp\ra$. By \pref{not2}(6), $R_{sh}=\rds+S_1\oplus\la S_2\ra$. Since $\la S_2\ra\sub \la R_{\lg}\ra$, it only remains to show that $\rds+S_1\sub\la\pp\ra.$ But as $\la S_1\ra=\sum_{i=1}^t\bbbz\sg_t$, we are done if we show that $\rds\sub\la\pp\ra$ and
$\sg_1,\ldots,\sg_t\in\la\pp\ra$.

By Lemma \ref{coron-100}, $\pp$ contains a short root $\a$. By Proposition \ref{modify1}(iv), we may assume that $\a\in\rds$. Since $\rdl\sub\la\pp\ra$, this implies that $\rds\sub\la\pp\ra$ and so $\dot R\sub\la\pp\ra.$ To conclude the proof for type $B_\ell$, it remains to
show that $\sg_1,\ldots,\sg_t\in\la\pp\ra$.
}

{For this, we fix $1\leq i\leq t$ and for $\a\in R^\times$ we consider the assignment
$$w_{\a}\stackrel{\psi_i}{\longmapsto}\left\{\begin{array}{ll}
	1&\hbox{if }\a\in\rds+\sg_i+\la L\ra,\\
	0&\hbox{otherwise}.
\end{array}\right.
$$
We claim that this assignment can be lifted to an epimorphism
$\w\rightarrow \bbbz_2$. By Theorem \ref{pbc1},
we need to show that the defining relations of the generalized presentation by conjugation 
vanish in $\bbbz_2$ under this assignment.
The relations of the form $w_\a w_\a$ are clear.
To check that $\psi_i(w_\a)\psi_i( w_\beta)\psi_i( w_\a)\psi_i( w_{w_\a(\beta)})=0$, we only need to versify that
$\beta\in\rds+\sg_i+\la L\ra$ if and only if $w_\a(\beta)\in\rds+\sg_i+\la L\ra$. 
So suppose 
$\beta=\dot\beta+\sg_i+\lam$ for some $\dot\beta\in\rds$
and $\lam\in\la L\ra$. Then 
$$w_\a(\beta)=w_\a(\dot\beta)+\sg_i+\lam \in
\left\{\begin{array}{ll}
	\pm\dot\beta+2\Lam+\sg_i+\lam\sub\rds+\sg_i+\la L\ra&\hbox{if }\a\in R_{sh},\\
	\rds+L+\sg_i+\lam\sub\rds+\sg_i+\la L\ra&\hbox{if }\a\in R_{lg},
\end{array}\right.
$$
as required.}

{Finally, we consider relations of the form  \pref{pgraph4}(\ref{(iii)}). For this, we must show that for 
	$$(\a,\sg=\sum_{j=1}^\nu m_j\sg_j)\in(\{\theta_s\}\times\Lam)\cup(\{\theta_\ell\}\times\sum_{i=t+1}^\nu\bbbz\sg_i),
	$$
	we have
	$$
	(\psi_i(w_{\a+\sg})\psi_i(w_\a))(\psi_i(w_{\a+\sg_1})\psi_i(w_\a))^{m_1}
	\cdots (\psi_i(w_{\a+\sg_\nu})\psi_i(w_\a))^{m_{\nu}}=0.
	$$
	Clearly, this holds if $\a=\theta_\ell$.
	Suppose now that $\a=\theta_s$. Since $\psi_i(w_{\theta_s})=0$ and
	$\psi_i(w_{\theta_s+\sg_j})=0$ for all $j\neq i$, we must show
	$\psi_i(w_{\theta_s+\sg})+m_i\psi_i(w_{\a+\sg_i})=0$. Now
	if $\sg\in\sg_i+\la L\ra$, then as $1\leq i\leq t$, we have $m_i\in2\bbbz+1$. Then $\psi_i(w_{\theta_s+\sg})=1$ and
	$m_i\psi_i(\theta_s+\sg_i)=1$ and we are done. Finally, we consider the case
	$\sg\not\in\sg_i+\la L\ra$. This forces $m_i\in2\bbbz$.
	Then we have 
	$\psi_i(w_{\theta_s+\sg})=0$ and
	$m_i\psi_i(\theta_s+\sg_i)=0$ which again gives the result.
	Thus as was claimed, the assignment $\psi_i$ induces
	an epimorphism $\psi_i:\w\rightarrow\bbbz_2$.}

{Now let
	$\a\in\rds$ and fix $1\leq i\leq t$. Then $\a+\sg_i\in R$ and so
	$w_{\a+\sg_i}=w_{\gamma_1}\cdots w_{\gamma_n}$ for some
	$\gamma_i\in\pp$. Then $\psi_i(w_{\a+\sg_i})=\psi_i(w_{\gamma_1}\cdots w_{\gamma_n}),$ and so 
	$$1=\psi_i(w_{\a+\sg_i})=\sum_{j=1}^{k}\psi_i(w_{\gamma_{i_j}}),$$
	where $\gamma_{i_j}\in\pp\cap R_{sh}$. Thus at least for one $j$, we get $\psi(w_{\gamma_{i_j}})=1$, implying that $\gamma_{i_j}\in(\rds+\sg_i+\la L\ra)\cap\pp$. Since $\la L\ra\sub\la\pp\ra$, this implies $\dot\gamma+\sg _i\in\la\pp\ra$ for some $\dot\gamma\in\rds$.  Since $\rds\sub\la\pp\ra$, we get $\sg_i\in\la\pp\ra$ and the proof for type $B_\ell$ is completed. The proof for type $C_\ell$ is analogous, replacing the roles of short and long roots.}

{Finally, we show that $\pp$ is connected. If not, then $\pp=\pp_1\uplus\pp_2$ where $\pp_1$ and $\pp_2$ are nonempty sets and $(\pp_1,\pp_2)=\{0\}$. Now let $\a\in R^\times$. Since $\pp$ spans $\V$, there exists $\beta\in\pp$ with $(\a,\beta)\not=0$. Assume without loss of generality that $\beta\in\pp_1$.
	We have $w_\a=w_{\a_1}\cdots w_{\a_m}$ for some
	$\a_1,\ldots,\a_m\in\pp$. Since $\pp_1$ and $\pp_2$ are orthogonal, it follows that
	$\beta-(\beta,\a^\vee)\a=w_\a(\beta)=w_{\a_{j_1}}\cdots w_{\a_{j_k}}(\beta)$ where $\a_{j_1},\ldots,\a_{j_k}$ are in $\pp_1$.
	Since $(\a,\beta)\not=0$, this forces $\a\in\pp_1$. This argument gives $R^\times\sub (R^\times\cap\span_{\bbbr}\pp_1)\uplus (R^\times\cap\span_{\bbbr}\pp_2)$ which is absurd as $R$ is connected.}\qed

To each subset $\pp$ of $R^\times$, we associate two subsets of $R$ as follows.

\begin{equation}\label{coron3}
R_\pp^{n}:=\w_\pp\pp\andd R_\pp^{i}:=(R_\pp^n-R_\pp^n)\cap R^0.
\end{equation}

\begin{pro}\label{coron2}
Let $\emptyset\neq\pp\sub R^\times$ be connected.

(i)
The set $R_\pp=R_\pp^n\cup R_\pp^i$ is an extended affine root system in
$\V_\pp:=\span_{\bbbr}R_\pp$ {satisfying} $R_\pp^\times=R_\pp^n$,             
$R_\pp^0=R_\pp^i$ and $\w_\pp=\w_{R_\pp}$. 

(ii) If $\w=\w_\pp$, then $R_\pp$ has the same rank, nullity and type of $R$.
In particular, if $\pp$ belongs to either of the classes
$\mm_c$ or $\mm_m$, then $R_\pp$ {has the same rank, nullity and type of $R$}. 

(iii) If $\pp$ is a reflectable set, then $R_\pp=R$. 


\end{pro}

\proof {First, we note from the way  $R_\pp$ is defined that $R^0_\pp=R^i_\pp$ and
$R^\times_\pp=R^n_\pp$.  

(i) From definition of $R_\pp$ it is clear that axioms (R1)-(R3) and (R5) of an extended affine root system (Definition \ref{def2}) hold for $R_\pp$. 
Since $\pp$ is connected it follows easily that  $R_\pp$ is connected, so (R6) also holds. For (R4), since $\V^0_\pp=\V_\pp\cap\V^0$, we have
\begin{eqnarray*}
	R^0_\pp=R^i_\pp&=&(R^n_\pp-R^n_\pp)\cap R^0\\
	&=&
	(R^\times_\pp-R^\times_\pp)\cap \big((R^\times-R^\times)\cap\V^0\big)\\
	&=&(R^\times_\pp-R^\times_\pp)\cap \V^0_\pp.
\end{eqnarray*}
}

(ii) {If $\w=\w_\pp$, then $\span_{\bbbr}\pp=\V$ and $\V_\pp=\V$, so
$\span_{\bbbr} R^\times_\pp=\span_{\bbbr}\pp=\V$ and
$\span_{\bbbr} R_\pp^0=\V^0$. Thus $R_\pp$ has the same nullity as $R$. This gives $\rank\; R_\pp=\dim \V_\pp-\dim \V^0_\pp=
\dim\V-\dim\V^0=\rank\;R$. Also by Lemma \ref{coron-100}, $R$ and $R_\pp$ have the same number of root lengths.}
{Thus both $\bar R$ and $\bar{R}_\pp$ have the same rank and the same number of root lengths, with $\bar{R}_\pp\sub\bar{R}$. Now consulting the root data of finite root systems (see \cite[Table 12.2.1]{Hu1}), we conclude that $\bar{R}_\pp=\bar{R}$.}

(iii) {Since $R_\pp^\times=R^n_\pp=R^\times$ and
	$R^0_\pp= (R^n_\pp- R^n_\pp)\cap R^0=(R^\times- R^\times)\cap \V^0=R^0$, we are done.}
\qed  

\subsection{The finite case}\label{the finite case}
We show that $\mm_r=\mm_c=\mm_m$ if $R$ is finite.

\begin{lem} \label{L1}
	Let $R$ be a reduced irreducible finite root system of rank $\ell$;
	then $c=\ell$.
\end{lem} 
\proof {Since the rank of a finite root system is equal to the minimal number of generators for the corresponding Weyl group, we get $\ell=c$.}

\begin{pro}\label{T3}
{Let $R$ be a reduced irreducible 
	finite root system, then $\mm_c=\mm_{m}=\mm_{r}$.}
\end{pro}

\proof
{Let {$\pp$} be a reflectable base for $R$. By Lemma \ref{0lem}, 
$|\pp|=\mbox{rank}(R)$. Thus by Lemma \ref{L1}, $\mm_r\subseteq\mm_c$.
Now let {$\pp\in\mm_c$; since $\w=\w_\pp$, we get from Proposition \ref{coron2}} that $R_\pp$ is a finite subsystem of $R$ of the same type and rank of $R$. Thus $R=R_\pp$
and so $\pp$ is a reflectable set for $R$. This proves that $\mm_c\sub\mm_r,$ and so $\mm_c=\mm_r$.}

Next, let $\pp\in\mm_m$. As we saw in the previous paragraph,
$\pp$ is a reflectable set. If $|\pp|>\rank\;R=c$, then $\pp$ contains a proper subset $\pp'$ which is a reflectable base, and so $\w_{\pp'}=\w$, contradicting that $\pp\in\mm_m$. Thus 
$$\mm_c\sub\mm_m\sub\mm_r\sub\mm_c,$$
completing the proof that $\mm_c=\mm_m=\mm_r.$\qed

The following example shows that Proposition \ref{T3} may fail if the extended affine root system $R$ is not finite.
\begin{exa}\label{ex1}
{Consider the extended affine root system  $R=\Lam\cup(\pm\ep+\Lam)$ of type $A_1$, where $\{0,\pm\ep\}$ is a finite root system of type $A_1$ and $\Lam=\bbbz\sg_1\oplus\bbbz\sg_2\oplus\bbbz\sg_3$. 
	Set
	$$\pp:=\{\ep, \sg_1-\ep, \sg_2-\ep, \sg_3-\ep, \sg_1+\sg_2-\ep, \sg_1+\sg_3-\ep, \sg_2+\sg_3-\ep,\sg_1+\sg_2+\sg_3-\ep\}.$$  
	By Proposition \ref{L11}(i), $\pp$ is a reflectable base for $R$, so $\pp\in\mm_{r}$. A simple computation shows that
	$$w_{\sg_1+\sg_2 +\sg_3-\ep}=w_{\sg_1+\sg_3-\ep}w_{\sg_1-\ep}w_{\sg_1+\sg_3-\ep}w_\ep w_{\sg_3-\ep}w_{\ep}w_{\sg_2-\ep}w_{\sg_2+\sg_3-\ep}w_{\ep},$$
	so $w_{\sg_1+\sg_2 +\sg_3-\ep}\in\la w_\a\mid \a\in\pp\setminus\{w_{\sg_1+\sg_2 +\sg_3-\ep}\}\ra$; thus $\pp\not\in\mm_m$ and then by (\ref{coron1}), $\pp\not\in\mm_c$. That is $\pp\not\in\mm_{rc}$.}
\end{exa}

\section{Reflectable bases, the class $\mm_r$}\setcounter{equation}{0}
In this section, we consider various aspects of reflectable bases; in particular, we investigate their existences and their possible cardinalities.

\subsection{Characterization theorem for reflectable bases}
We recall two recognition theorems for reflectable bases associated with reduced extended affine root systems given in \cite[Theorems 3.1, 3.14, 3.24, 3.26, 3.27]{AYY} and \cite[Theorem 3.2]{ASTY1}.

\begin{thm}\label{coron101} (Recognition Theorem for simply laced types) Let $R$ be  simply laced of type $X$ and  $\Pi\sub R^\times$ satisfy $\la\Pi\ra=\la R\ra.$ Then $\Pi$ is a reflectable  base for $R$ if and only if
	
	(i) $R^\times=\uplus_{\a\in\Pi}\big((\a+2\la R\ra)\cap R\big)$, if $X=A_1$,
	
	(ii) $\Pi$ is a minimal generating set for the free abelian group $\la R\ra$, if $X$ is simply laced of rank$>1$.
\end{thm}

\begin{thm}\label{coron101-1} (Recognition Theorem for non-simply laced types) Let $R$ be non-simply laced of type $X$ and
	{$\Pi\sub R^\times$ satisfy $\la \Pi\ra=\la R\ra$.} Then $\pp$ is a reflectable base for $R$ if and only if
	
	(i) $R _{sh}=\uplus_{\a\in\Pi_{sh}}\big((\a+\la R _{lg}\ra)\cap R _{sh}\big)$, and
	$R _{lg}=\uplus_{\a\in\Pi_{lg}}\big((\a+ 2\la R _{sh}\ra)\cap R _{lg}\big),$ if $X=B_2$,
	
	(ii) $R_{sh}=\uplus_{\a\in\Pi_{sh}}\big((\a+\la R_{lg}\ra)\cap R_{sh}\big)$, and $\Pi_{lg}$ is a minimal set with respect to the property that  $\{\a +2\la R_{sh}\ra\mid \a\in \Pi_{lg}\}$  is a basis for the $\bbbz_2$-vector space $\la R_{lg}\ra/2\la R_{sh}\ra$, if $X=B_\ell,$ $\ell\geq 3$,
	
	(iii) $R_{lg}=\uplus_{\a\in\Pi_{lg}}\big((\a+2\la R_{sh}\ra)\cap R_{lg}\big)$, and  $\Pi_{sh}$ is  minimal with respect to the property that the set $\{\a+\la R_{lg}\ra\mid\a\in\Pi_{sh}\}$ is a basis for the $\bbbz_2$-vector space $\la R_{sh}\ra/\la R_{lg}\ra,$
	if $X=C_\ell$, $\ell\geq 3$,
	
	(iv) 
	$\Pi_{sh}$ is minimal with respect to the property that the set $\{\a +\la R_{lg}\ra\mid\a\in\Pi_{sh}\}$ is a basis for the $\bbbz_2$-vector space $\la R_{sh}\ra/\la R_{lg}\ra$, and
	$\Pi_{lg} $ is minimal with respect to the property that the set $\{\a +k\la R_{sh}\ra\mid\a\in\Pi_{lg}\}$ is a basis for the $\bbbz_2$-vector space $\la R_{lg}\ra/k\la R_{sh}\ra$, if $X=F_4$, $G_2$, where $k$ is as in \pref{not2}(4).
\end{thm}

\begin{rem}\label{coron33}
	If in Theorems \ref{coron101} and Theorem \ref{coron101-1}, we drop the terms ``minimal'' and ``basis'' and change the symbol ``$\uplus$'' to ``$\cup$'', then the conditions there characterize the reflectable sets in $R$.
\end{rem}

\subsection{Concrete family of reflectable bases}
We show that set $\pp(X)$ given in Table \ref{tab:fq} is a reflectable base for the extended affine root system of type $X$.

\begin{pro}\label{L11}
	Let $R$ be a reduced extended affine root system of type $X$, nullity $\nu$, and twist number $t$. Then the following hold. 
	
	(i) the set $\pp(X)$ given in Table \ref{tab:fq} is a reflectable base for $R$ with $|\pp(X)|=\mar{ind}(R)+\ell+\nu$. 
	
	(ii) If $\ind(R)=0$, then $\pp(X)\in\mm_c$.
\end{pro}
\proof
{We follow the same notations as in \pref{not2}. 
	
	(i) Suppose first that $R$ is simply laced of rank $>1$, or is of type $F_4$ or $G_2$. Since
	$|\pp(X)=\ell+\nu$, we get from
	Theorem \ref{coron101}(ii) and Theorem \ref{coron101-1}(iv) that $\pp(X)$ is a reflectable base.

	Next suppose $X=B_\ell$ with $\ell\geq 3$. We have 
	$$\pp:=\pp(B_\ell)=\{\a_1,\ldots, \a_\ell, \gamma_1-\te_s, \ldots, \gamma_{\rho_1}-\te_s, \d_{t+1}-\te_l, \ldots, \d_\nu -\te_l\},$$
	where $\rho_1=\ind(S_1)$. 
	We assume that the only short root of $\dot\Pi$ is $\a_\ell$. 
	We have 
	$$\begin{array}{l}
		\pp_{sh}=\{\a_\ell, \gamma_1-\te_s, \ldots, \gamma_{\rho_1}-\te_s\},\\
		\pp_{lg}=\{\a_1, \ldots, \a_{\ell -1}, \d_{t+1}-\te_l, \ldots, \d_\nu -\te_l \}.
	\end{array}
	$$ 
	By Proposition \ref{L11}, $\pp$ is a reflectable set for $R$ and so by Remark \ref{coron33}, 
	
	(1) $R_{sh}=\cup_{\a\in\pp_{sh}}(\a+\la R_{lg}\ra)\cap R_{sh}$,
	
	(2) $\{\a +2\la R_{sh}\ra\mid \a\in \pp_{lg}\}$ {spans} the $\bbbz_2$-vector space $\u={\la R_{lg}\ra}/{2\la R_{sh}\ra}$.}

\noindent Thus to show that $\pp$ is a reflectable base, it is enough to show (by Theorem \ref{coron101-1}(ii)) that the union in (1) is disjoint and that the set in (2) is a basis for $\u$.  We start by proving that the union in (1) is disjoint. 

If $\a_\ell\in\cup_{\a\in\pp_{sh}\setminus\{\a_\ell\}}(\a+\la R_{lg}\ra)\cap R_{sh}$, then there is an $1\leq i\leq \mbox{ind}(S_1)$ such that 
$\a_\ell\in\gamma_i -\te_s +\la R_{lg}\ra$ and this implies that $\gamma_i\in\la L\ra$, which contradicts \pref{not2}(\ref{s1}).
Next suppose there exists an $1\leq i\leq \mbox{ind}(S_1)$ such that
$\gamma_i -\te_s\in\cup_{\a\in\pp_{sh}\setminus\{\gamma_i -\te_s\}}(\a+\la R_{lg}\ra)\cap R_{sh}$.
Now if $\gamma_i -\te_s\in \a_\ell +\la R_{lg}\ra$, then $\gamma_i\in\la L\ra$ which again contradicts \pref{not2}(\ref{s1}). If $\gamma_i -\te_s\in\gamma_j -\te_s +\la R_{lg}\ra$ for some $1\leq j\neq  i\leq\mbox{ind}(S_1)$, we get $\gamma_i -\gamma_j\in 2\la S_1\ra$, which contradicts \pref{not2}(\ref{coron6}). {This completes the proof for (1).} 

Next, we consider (2). {As  the set
	$\{\a +2\la R_{sh}\ra\mid \a\in \pp_{lg}\}$ spans $\u$ and}
$\mbox{dim}_{\bbbz_2}\u=(\ell -1)+(\nu -t)=|\pp_{lg}|=(\ell -1)+(\nu -t)$, {this set a basis for $\u$. Thus $\pp$ is a reflectable base for $R$  and $|\pp|=\mbox{ind}(R)+\ell +\nu$.}

{The proof for type $C_\ell(\ell\geq 3)$ is analogous to type $B_\ell$, replacing the roles of short and long roots in the proof, and using the characterization of reflectable bases of type $C_\ell$ given in Theorem \ref{coron101-1}(iii).}
The proofs for types $A_1$ and $B_2$ can also be carried out in a similar manner, using Recognition Theorems \ref{coron101}(i) and \ref{coron101-1}(i). 

{(ii)  Let $\pp=\pp(X)$. If $\mbox{ind}(R)=0$, then $|\pp|=\ell +\nu$ and as $\w=\w_\pp$, we have $|\pp|=\ell+\nu=\dim(\V)=\dim(\mbox{span}_\bbbr\pp)$; thus $\pp$ is a basis for $\V$. Since $\mbox{dim}(\V)\leq c\leq |\pp|$, we have $c=|\pp|=\ell +\nu$, that is, $\pp\in\mm_{c}$.}
\qed

For the proof of the following corollary, we recall this elementary fact from group theory that  if a group $G$ is finitely generated, then any generating set for $G$ contains a finite generating set.    

\begin{cor}\label{T1} 
{Let $\pp$ be in either of $\mm_c$, $\mm_m$ or $\mm_r$. Then
$\pp$ is finite.
}
\end{cor}
\proof {By \pref{pgraph9}(\ref{coron1}), we only need to consider $\pp\in\mm_m$ or  $\pp\in\mm_r$. By Proposition \ref{L11}, $\w$ is generated by reflections based on the finite set $\pp(X)$, and so is finitely generated. Thus any generating set for $\w$ contains a finite generating set. In particular, if $\pp\in\mm_m$ then $\pp$ contains a finite subset $\pp'$ such that $S_{\pp'}$ generates $\w$. But the minimality of $\pp$ gives, $\pp=\pp'$, and we are done.}
%

{Next, assume $\pp\in\mm_r$. We have $\w=\w_\pp$, and so as explained above, $\w$ is generated by reflections based on a finite subset $\pp'$ of $\pp$. Again let $\pp(X)$ be the reflectable base given in Proposition \ref{L11}. Then
$\pp(X)\sub R^\times=\w_{\pp}\pp=\w_{\pp'}\pp$. As $\pp(X)$ is finite, we have $\pp(X)\sub\w_{\pp'}\pp''$ for some finite subset $\pp''$ of $\pp$. Now setting {$\widehat\pp:=\pp'\cup\pp''$}, we get
$$R^\times=\w\pp(X)=\w_{\widehat\pp}\pp(X)\sub\w_{\widehat\pp}\w_{\pp'}\pp''\sub\w_{\widehat\pp}\widehat\pp.$$
Now the minimality of $\pp$ gives $\widehat\pp=\pp$, and we are done.}
\qed

\subsection{On cardinality of reflectable bases; types
	$A_1$, $B_\ell$, $C_\ell$, $F_4$, $G_2$}
In Corollary \ref{T1}, we just saw that reflectable bases have finite cardinalities.
If any two reflectable bases for an extended affine root system have {the same cardinality?} As we will see in the sequel the answer is not positive in general, but we will discuss situations in which the response is affirmative.  
\begin{thm}\label{T2} 
	Let $R$ be an extended affine root system of rank $\ell$, nullity $\nu$, and twist number $t$, and let $\pp$ and $\pp'$ be two reflectable bases for $R$. Then for the types given in Table \ref{tab:fq2} below, $|\pp|=|\pp'|$. Moreover, if $\pp\in\mm_r$, then $|\pp|$, $|\pp_{sh}|$ and $|\pp_{lg}|$ are given by Table  \ref{tab:fq2}, {where $S$, $S_1$ and $S_2$ are as in \pref{note2}(\ref{s1}).}
	
	\begin{table}[ht]
		\caption{Cardinality of reflectable bases}
		\vspace{-4mm}
		\label{tab:fq2}
		\begin{tabular}
			[c]{|c|c |c|c|} 
			\hline
			$\mbox{Type}$& $|\pp|$ & $|\pp_{sh}|$ & $|\pp_{lg}|$ \\
			\whline
			$A _1$ & $1+ \mbox{ind}(S)$ & --- & --- \\
			\hline
			$B_2$ & $2+ \mar{ind}(S_1)+\mar{ind}(S_2)$ & $1+ \mar{ind}(S_1)$ & $1+\mar{ind}(S_2)$\\
			\hline 
			$B_\ell$ $(\ell\geq 3)$ & $\ell+\mar{ind}(S_1)+(\nu -t)$ & $1+ \mar{ind}(S_1)$ & $(\ell- 1)+(\nu- t)$\\
			\hline
			$C_\ell$ $\ell\geq 3$ & $\ell +t+\mar{ind}(S_2)$ &$(\ell- 1)+ t$ & $1+ \mar{ind}(S_2)$\\
			\hline
			$F_4$ & $4+ \nu$ & $2+ t$ & $2+ (\nu- t)$ \\
			\hline
			$G_2$ & $2+ \nu$ & $1+ t$ & $1+ (\nu- t)$\\
			\hline
		\end{tabular}
		\medskip
	\end{table}
\end{thm}

\proof { If $R$ is of type $A_1$, then we have $|\pp|=1+\ind(S)$, by \cite[Lemma 1.7]{AP}, and we are done in this case.}

{Next suppose $R$ is of type $B_2$ and $\pp$ and $\pp'$ are two reflectable bases for $R$. By Theorem \ref{coron101-1}(i), we have the following coset description of $R_{sh}$ and $R_{lg}$:
$$R _{sh}=\uplus_{\a\in\pp_{sh}}(\a+\la R _{lg}\ra)\cap R _{sh}= \uplus_{\a'\in\pp'_{sh}}(\a'+\la R _{lg}\ra)\cap R _{sh}$$
and
$$R _{lg}=\uplus_{\a\in\pp_{lg}}(\a+ 2\la R _{sh}\ra)\cap R _{lg}= \uplus_{\a'\in\pp'_{lg}}(\a'+ 2\la R _{sh}\ra)\cap R _{lg}.$$
It follows that $|\pp_{sh}|=|\pp'_{sh}|$ and $|\pp_{lg}|=|\pp'_{\lg}|$. Thus $|\pp|=|\pp'|.$}

{Assume next that $R$ is of type $B_\ell$ $(\ell\geq 3)$ or  $C_\ell$ $(\ell\geq 3)$. We give the proof for type $B_\ell$; the proof for type $C_\ell$ is analogous, replacing the roles of short and long roots. Assume $\pp$ and $\pp'$ are two reflectable bases for $R$. According to Theorem \ref{coron101-1}(ii),  we have 
$$R _{sh}=\uplus_{\a\in\pp_{sh}}(\a +\la R _{lg}\ra)\cap R _{sh}=\uplus_{\a'\in\pp'_{sh}}(\a' +\la R _{lg}\ra)\cap R _{sh}$$
and
$$|\pp_{lg}|=\dim_{\bbbz_2}
\big({\la R_{lg}\ra}/{2\la R_{sh}\ra}\big)=|\pp'_{lg}|.$$
Thus $|\pp_{sh}|=|\pp'_{sh}|$ and $|\pp_{lg}|=|\pp'_{lg}|$, as required.} 

{Next assume $R$ is of type $F_4$ or $G_2$. By Theorem \ref{coron101-1}(iv), if $\pp$ and $\pp'$ are two reflectable bases for $R$,
then
$$|\pp_{lg}|=\dim_{\bbbz_2}
\big({\la R_{lg}\ra}/{2\la R_{sh}\ra}\big)=|\pp'_{lg}|\andd
|\pp_{sh}|=\dim_{\bbbz_2}
\big({\la R_{sh}\ra}/{\la R_{lg}\ra}\big)=|\pp'_{sh}|.$$}

Finally, we consider the last assertion in the statement. Let $\pp$ be a reflectable base of type $X=A_1,B_\ell, C_\ell$, $F_4$ or $G_2$.
Since now we know that any two reflectable bases for $R$ have the same cardinality, it is enough to assume that $\pp=\pp(X)$, where $\pp(X)$ is the reflectable base given in Proposition \ref{L11}. The result now can be seen from the information given in Table \ref{tab:fq}.
\qed 

\begin{cor}
	\label{cornew}
The	cardinality of a reflectable base is an isomorphism invariant of an extended affine root system of types  $X=A_1$, $B_\ell$, $C_\ell$, $F_4$, $G_2$.   
\end{cor}

\proof One notes that the index of semilattices $S$, $S_1$, and $S_2$ appearing in the structure of an extended affine root system, as well as the nullity, the rank, and the twist number are isomorphism invariants for an extended affine root system, see \cite[Chapter II]{AABGP}. Now the result is immediate from Proposition \ref{T2}.\qed

 \begin{lem} \label{min ref.}
 {Suppose $\pp$ is a reflectable base for $R$ such that $|\pp|=\ell +\nu$, then $\pp\in\mm_{rm}$ (that is $\pp$ is a minimal reflectable base).}	
 \end{lem}
\proof 
Since $\pp$ is a reflectable base, the reflections based on $\pp$ generate $\w$, so it is enough to show that no proper subset of $\pp$ has his property. But this is clear as $|\pp|=\dim \V$.
\qed

\subsection{On cardinalities of reflectable bases; simply laced types of rank $>1$}

{Let $R$ be a simply laced extended affine root system of rank $\ell>1$ and nullity $\nu>1$. We start with the following example which shows that Proposition \ref{T2} does not hold in general for $R$, namely reflectable bases for $R$ might have different cardinalities; even when
$\nu=1$ which is the affine case. The example also shows that  $c=\ell +\nu$ and that $\mm_{rc}$ is a proper subclass of $\mm_{rm}$. We then show that $\mm_r=\mm_{rm}$. The subsection is concluded with a technical result that shows that each reflectable base for $R$ contains a subset of cardinality $\ell+\nu$ which is a reflectable base for a subsystem of $R$ of the same type, rank, and nullity as $R$. }
\begin{exa}\label{ex2}	
{Let $R$ be a simply laced extended affine root system of type $X$ and rank$>1$. Then $R=\Lam\cup (\rd+\Lam)$, where $\rd$ is a finite root system of type $X$ and $\Lam=\bbbz\sg_1\oplus\cdots\oplus\bbbz\sg_\nu$, see \pref{pgraph2} and \pref{not2}. Let $\dot{\Pi}=\{\a_1,\ldots,\a_\ell\}$ be a fundamental system for $\rd$.
	Set	$$\Pi=\{\a_1,\ldots,\a_\ell,\a_1+m_1\sg_1,\a_2+m_2\sg_1, \a_2+\sg_2,\ldots,\a_2+\sg_\nu\},$$
	where $m_1,m_2>1$ are relatively prime integers. 
	We claim that $\Pi$ is a reflectable base for $R$. Since $(m_1,m_2)=1$, we have $\la\Pi\ra=\la\dot R\ra\oplus\Lam=\la R\ra$, and so we get from Theorem \ref{coron101}(ii) and Remark \ref{coron33}, that $\Pi$ is a reflectable set for $R$. Therefore, using again Theorem \ref{coron101}(ii), it is enough to show that for each $\a\in\Pi$, $\Pi'_\a:=\Pi\setminus\{\a\}$ does not generate $\la R\ra$. Clearly, we only need to check this for $\a\in\{\a_1,\a_2,\a_1+m_1\sg_1,\a_2+m_2\sg_1\}$. Since $m_1,m_2>1$, we have $\sg_1\not\in \la\Pi'_{\a_1+m_1\sg_1}\ra$ and similarly $\sg_1\not\in \la\Pi'_{\a_2+m_2\sg_1}\ra$. If $\sg_1\in\Pi'_{\a_1}$, then $\sg_1=k\a_2+k'(\a_1+m_1\sg_1)+k''(\a_2+m_2\sg_1)$ for some $k,k',k''\in\bbbz$. This gives $k'=0$ and so $k''m_2=1$ which is impossible. The same reasoning gives $\sg_1\not\in\la\Pi'_{\a_2}\ra.$
	Thus $\Pi$ is a reflectable base. Note that $|\Pi|=\ell+\nu+1>|\pp(X)|=\ell+\nu$, where $\pp(X)$ is the reflectable base given in Table \ref{tab:fq}.} We conclude that $\mm_{rc}$ is a proper subclass of
$\mm_{rm}$.
\end{exa}

\begin{pro}\label{min AD}
	{Let $R$ be a simply laced extended affine root system of rank $>1$. Then $\mm_r=\mm_{rm}$.}		
\end{pro}

\proof {We must show $\mm_r\sub\mm_m$. Assume that $\pp$ is a reflectable base for $R$; by Theorem \ref{coron101}(ii), $\pp$ is  minimal with respect to the property that $\la\pp\ra=\la R\ra.$
	If $\pp\not\in\mm_m$, then there is at least one $\a_0\in\pp$ such that $w_{\a_0}\in\la w_\ga |\ga\in\pp\setminus\{\a_0\}\ra$. 
	But then Proposition \ref{coron-21}, gives $\la\pp\setminus\{\a_0\}\ra=\la R\ra$, which contradicts the minimality of $\pp$.}\qed

\begin{pro}\label{pro A}
	Let $R$ be a simply laced extended affine root system of rank $\ell>1$ and nullity $\nu$. Let $\Pi$ be a reflectable base for $R$. Then $\Pi$ contains a subset $\Pi'$ of cardinality $\ell+\nu$ such that $\Pi'$ is a reflectable base for a root subsystem $R'$ of $R$, where $R'$ has the same type, rank, and nullity as $R$.
\end{pro} 

\proof Let $\Pi$ be a reflectable base for $R$, then $|\Pi|\geq\ell +\nu$. If $|\Pi|=\ell +\nu$, there is nothing to prove. Now assume that $|\Pi|> \ell +\nu$. Recall from Subsection \ref{extended affine root systems} that $\bar{\;}$ is the canonical map from $\V$ onto $\bar{\V}:=\V\backslash\V^0$ and $\bar{R}$ is the image of $R$ under $\bar{\;}$. Since {$\w_{\Pi}\Pi=\rcross$}, we get from Proposition \ref{modify1} that there exists a subset  $\dot\Pi$ of $\Pi$, such that
$\dot R=(\w_{\dot\Pi}\dot\Pi)\cup\{0\}\sub R$ is an irreducible  finite root system isomorphic to $\bar R$, and
we get the description $R=\Lam\cup(\dot R+\Lam)$ for $R$ where $\Lam$ is a lattice, see \pref{pgraph2}. 
We fix a reflectable base $\{\a_1,\ldots,\a_\ell\}\sub\dot\Pi$ of
$\dot R$, see \pref{the finite case}. Then
$\dot\V=\span_{\bbbr}\dot R$. We extend $\{\a_1,\ldots,\a_\ell\}$ to a basis $\Pi'=\{\a_1,\ldots,\a_\ell,\a_{\ell+1},\ldots,\a_{\ell+\nu}\}$ of $\V$, such that $\Pi'\sub\Pi$.  Since $\Pi'$ spans $\V$, it is connected.

We now consider the extended affine root system $R':=R_{\Pi'}$ defined
in Proposition \ref{coron2}. {By the same argument as in the proof of Proposition \ref{coron2}(ii)}, $R'$ has the same type, rank and nullity as $R$. Since $|\Pi'|=\ell+\nu=\dim \V$, it is a reflectable base for $R'$.\qed

 \section{A characterization of minimal root systems}\label{Minimal bases} \setcounter{equation}{0}
As in the previous sections, we let  $R$ be an extended affine root system of reduced type $X$ with corresponding Weyl group $\w$. It is known that $\w$ has the presentation by conjugation, see Definition \ref{pbc}, if and only if $R$ is "minimal'' in the sense of Definition \ref{R-min} below; therefore minimality condition reflects the geometric aspects of extended affine root systems. In this section, we characterize minimal root systems in terms of minimal reflectable bases.

\subsection{Minimal extended affine root systems}
\begin{DEF}\label{R-min}
{Following \cite{Ho}, {we say the extended affine root system $R$ is {\it minimal},} if for each $\a\in R^\times$, 
$\w_\pp\subsetneq\w$, where $\pp=R^\times\setminus\w\a$.}
\end{DEF}

\begin{lem}\label{equiv1}
Let $\pp\sub R^\times$.	 	

{(i)  $\w\pp=\pp$  if and only if $S_\pp=S_{\w\pp}$ and $\pp=-\pp$.}
%
 	
{(ii) The extended affine root system $R$ is  minimal if and only if whenever $\pp\sub R^\times$ with 
$\w_\pp=\w$ and $\w\pp=\pp$, then $\pp=R^\times$.} 
%
\end{lem}
\proof {(i) The ``if'' part is clear.
{Now assume that $S_\pp=S_{\w\pp}$ and $\pp=-\pp$. Then for $w'\in\w$ and $\a\in \pp$ we have $w_{w'(\a)}\in S_\pp$; thus there exists  $\be\in \pp$ such that $w_{w'(\a)}=w_\be$ and then we have $w'(\a)=\pm\be$. Now as $\pp=-\pp$, we have $w'(\a)\in \pp$ and therefore $\w \pp=\pp$.}  

{(ii) Assume first that  $R$ is  minimal and $\pp$ is a subset of $\rcross$ such that $\w=\w_\pp$ and $\w\pp=\pp$. We must show $\pp=R^\times.$ If not, we pick $\a\in R^\times\setminus\pp$. Since $\w\pp=\pp$, we get $ \pp\subset R^\times\setminus\w\a$. But now the assumption $\w=\w_\pp$ contradicts minimality of $R$.}

{Next assume that the ``only if'' part holds, and that $R$ is not minimal. So there exists $\a\in R^\times$ such that $\w_\pp=\w$ for  $\pp:=R^\times\setminus\w\a$. Now since $\w\pp=\pp$, we get by the assumption that $R^\times=\pp$, which is absurd. }
\qed

\begin{thm}\label{L11-1}
	Let $R$ be a reduced extended affine root system of type $X$ and $\pp(X)$ be as in Table \ref{tab:fq}. Then  $\pp(X)\in\mm_m$ if and only if $R$ is a minimal root system.
\end{thm}

\proof
First, suppose that $R$ is of type $A_1$ or one of the non-simply laced types; then by \cite[Theorem 5.16]{AS3} and \cite[Theorem 2.33]{AS5}, $\pp(X) \in\mm_{m}$ if and only if $R$ is a minimal root system. {Next, assume that $R$ is simply laced of rank $>1$ and $\pp\in\mm_{m}$. If $R$ is not minimal, then by \cite[Theorem 4.5]{Ho}, $\w$ does not have the presentation by conjugation, contradicting Theorem \ref{pbc1}(iii). 
	
	Conversely, assume that $R$ is a minimal root system. As $|\pp(X)|=\ell+\nu=\mbox{dim}(\V)$, we have $\pp(X)\in\mm_{c}\sub\mm_{m}$.\qed
	
	\subsection{Minimality and presentation by conjugation}
	Let $R$ be a reduced extended affine root system with Weyl group $\w$. Let $\widehat\w$ be the group defined by 
	generators $\widehat{w}_\a$, $\a\in R^\times$ and relations 
	
	- $\widehat{w}^2=1$, $\a\in R^\times$,
	
	- $\widehat{w}_\a \widehat{w}_\beta \widehat{w}_\a=\widehat{w}_{w_\a\beta}$, $\a,\beta\in R^\times$.
	
	We recall from Definition \ref{pbc} that if $\w\cong\widehat\w$, then $\w$ is said to have the presentation by conjugation. We refer to this group as the {\it conjugation presented group}, associated with $R$.

	\begin{pro}\label{min-gen.} Given any reflectable base $\pp$ for $R$,
		the set $\{\widehat{w}_\a\mid \a\in\pp\}$ is a minimal generating set for the conjugation presented group $\widehat\w$ associated to $R$.
	\end{pro}
	\proof
	{If $\nu=0$ then $R$ is a finite root system and so $\w\cong\widehat\w$. Then we are done by  Theorem \ref{T3}.} 
	
	{Assume next that $\nu\geq 1$.  If $R$ is simply laced  of rank$>1$ or has one of types $F_4$ or $G_2$, then by Lemma \ref{min ref.}, Propositions \ref{T2} and \ref{min AD},  $\{w_\a\mid\a\in\pp\}$ is a minimal generating set for $\w$. But by Theorem \ref{pbc1}, for the types under consideration $\w\cong\widehat\w$ and so we are done.}
	
	{Next, suppose that $R$ is of type $A_1$; then by Theorem \ref{coron101}(i), we have
		\begin{equation}\label{temp1}
			\rcross=\uplus_{\a\in\pp}(\a +2\la R\ra)\cap\rcross.
		\end{equation}
		On the other hand, since $R^\times=\cup_{\a\in\pp}\w\a$, and for $\a\in R^\times$ we have $\w\a\subseteq \a+ 2\la R\ra$, we get from
		(\ref{temp1}),
		\begin{equation}\label{temp2}
			\rcross=\uplus_{\a\in\pp}\w\a.
		\end{equation}  
		If $\{\widehat{w}_\a\mid\a\in\pp\}$ is not a minimal generating set for $\widehat{\w}$; then there is at least one 
		$\be\in\pp$ such that $\widehat{w}_{\be}=\widehat{w}_{\ga_1}\cdots \widehat{w}_{\ga_n}$ for some $\ga_i\in\pp\setminus\{\be\}$, $1\leq i\leq n$. By Corollary \ref{cor12}, the assignment
		$$
		w_\a\stackrel{\Phi_\be}{\longmapsto}\left\{\begin{array}{ll}
			1 & \mbox{if}~ \a\in\w\be,\\
			0 & \mbox{if} ~\a\in\rcross\setminus\w\be,
		\end{array}\right.
		$$
		induces an epimorphism $\Phi_\b:\w\rightarrow\bbbz_2$.
		Now, by (\ref{temp2}), $\pp\setminus\{\be\}\subseteq \rcross\setminus\w\be$, and so we get ${\Phi}_\be(\widehat{w}_\be)=1$ and 
		${\Phi}_\be(\widehat{w}_{\ga_1}\ldots \widehat{w}_{\ga_n})=0$
		which is a contradiction.}

{Next, we consider type $B_2$. In this case, it is easy to see that
	\begin{equation}\label{temp3}
		\w\a\subseteq \left\{\begin{array}{ll}\a +\la R_{lg}\ra &\hbox{if }\a\in R_{sh},\\ 
			\a +2\la R_{sh}\ra&\hbox{if }\a\in R_{lg}.\\
		\end{array}\right.
	\end{equation}
	By Theorem \ref{coron101}(iii), we have
	$$R _{sh}=\uplus_{\a\in\pp_{sh}}\big((\a+\la R _{lg}\ra)\cap R _{sh}\big)\andd
	R _{lg}=\uplus_{\a\in\pp_{lg}}\big((\a+ 2\la R _{sh}\ra)\cap R _{lg}\big).
	$$
	Using these equations together with (\ref{temp3}), we get that the union
	$R^\times=\cup_{\a\in\pp}\w\a$ is disjoint. Then using the same argument as in type $A_1$, we achieve the desired result.}

{Finally, suppose that $R$ is of type $B_\ell$, $(\ell>2)$, and assume to the contrary that $\{\widehat{w}_\a|\a\in\pp\}$ is not a minimal generating set for $\widehat{\w}$; then there is at least one $\be\in\pp$ such that $\widehat{w}_\be\in\la\widehat{w}_\a \mid\a\in\pp\setminus\{\be\}\ra$. We proceed with the proof by considering the following two cases.}

{$\be\in R_{sh}$:  Similar to the case $B_2$, using Theorem \ref{coron101}(iv), we see that $R_{sh}=\uplus_{\a\in\pp_{sh}}\w\a$, then $\pp\setminus\{\beta\}\sub R^\times\setminus\w\be$.  Therefore repeating the same argument as in type $B_2$ leads to a contradiction.}

{ $\be\in R_{lg}$: We have $\widehat{w}_\be=\widehat{w}_{\ga_1}\cdots \widehat{w}_{\ga_n}$ for some $\gamma_i\in\pp\setminus\{\be\}$.  Applying  the epimorphism $\psi$ given in \ref{pgraph7}(\ref{epi})), we get $w_\be=w_{\ga_1}\ldots w_{\ga_n}$. Since $\ell\geq 3$, there exists $\a\in\rlg$  such that $(\a, \be^\vee)=-1$. If $\{\ga_{i_1}, \ldots\ga_{i_m}\}=\{\ga_1, \ldots, \ga_n\}\cap\rlg$ then we have}
$$
\a +\be=w_\be(\a)=w_{\ga_1}\ldots w_{\ga_n}(\a)
\in w_{\ga_{i_1}}\ldots w_{\ga_{i_m}}(\a)+2\la \rsh\ra;
$$
thus 
\begin{equation}\label{equ2}
	\be\in \la\ga_{i_1}, \ldots ,\ga_{i_m}\ra +2\la\rsh\ra,\quad\gamma_{i_j}\in\pp_{lg}\setminus\{\be\},
\end{equation}
which is a contradiction as by Theorem \ref{coron101}(iv), $\pp_{lg}$ is minimal with respect to the property that  the set $\{\a+ 2\la\rsh\ra\mid\a\in\pp_{lg}\}$ is a basis for the $\bbbz_2$ vector space ${\la\rlg\ra}/{2\la\rsh\ra}$.

{As usual, the argument for the case $C_\ell$ is analogous to $B_\ell$, replacing the roles of short and long roots.}
\qed

\begin{thm}\label{pro13} {The following are equivalent:}
	
	(i) $R$ is minimal,
	
	(ii) $\w$ has the presentation by conjugation,
	
	(iii)  $\mm_r=\mm_{rm}$.
\end{thm}

\proof Conditions (i) and (ii) are equivalent by \cite[Theorem 5.8]{Ho}. 
Next assume (ii) holds, namely $\w\cong\widehat\w$ under the isomorphism induced by  the assignment $\psi:\widehat\w\rightarrow\w$, $\widehat{w}_\a\mapsto w_\a$, $\a\in R^\times$. 
Let $\Pi\in\mm_r$.  By Proposition \ref{min-gen.}, the set $\{\widehat{w}_\a\mid\a\in \Pi\}$ is a minimal generating set for $\widehat\w$, and so we are done. 
This gives (iii).

Finally, assume that (iii) holds.  Consider the set $\pp(X)$ given in Table \ref{tab:fq}. By Proposition \ref{L11}, $\pp(X)\in\mm_r$ and so $\pp(X)\in\mm_m$. Then $R$ is minimal by {Theorem \ref{L11-1}}. This gives (i) and the {proof is completed}.\qed

\begin{cor}\label{cor13} Let $R$ be a reduced extended affine root system of 
type $X$, rank $\ell$ and nullity $\nu$. {Then under either of the following two conditions:} 

(i) {$\nu\leq 2$,}

(ii) {$X=A_\ell$, $(\ell>1),$ $D_\ell$, $E_{6,7,8}$, $F_4$ or $G_2$},

\noindent {\it we have $\mm_r=\mm_{rm}$}. Moreover, if $\ind(R)=0$ and $X=A_1$, $B_\ell,\; (\ell\geq 2)$ or $C_\ell,\; (\ell\geq 3)$, then $\mm_r=\mm_{rc}=\mm_{rm}$.
\end{cor}
\proof For the first part of the statement, we only need to show by Theorem \ref{pro13} that under either of the conditions (i) and (ii), $\w$ has the presentation by conjugation. But this is immediate from Theorem \ref{pbc1}(ii-(iii). 

Next, we assume $\ind(R)=0$ where
$X=A_1$, $B_\ell$, $C_\ell$. By \pref{pgraph9}(\ref{coron1}), it is enough to show that $\mm_r\sub\mm_{rc}$. Let $\pp\in\mm_r$. 
By Proposition \ref{T2}, we have $|\pp|=\ell +\nu$. Now as $\w=\w_\pp$, we get $\V=\span_{\bbbr}\pp$, and since by \pref{pgraph9}(\ref{coron11}), $\dim(\V)\leq c\leq |\pp|$, we get $c=\ell +\nu$. Thus $\pp\in\mm_{rc}$.  \qed

\subsection{Gallery of results}
{The following table, Table 3, concludes and illustrates our investigation concerning the interrelations between classes $\mm_r,\mm_m$ and $\mm_c$.} Here $R$ is an irreducible reduced extended affine root system of type $X$, rank $\ell$ and  nullity $\nu$.

{\scriptsize
	\begin{table}[ht]\label{test}\caption{Gallery of results  on the level of root system and Weyl group}
		\vspace{-4mm}
{\begin{tabular}{|c|c|c|}
\hline
$R=\hbox{ extended affine root system}$ &$\stackrel{\hbox{Proposition \ref{L11}}}{{\Longrightarrow}}$ & $\pp(X)$ is reflectable base for $R$ \\
\hline
$\pp(X)\in\mm_m$ &$\stackrel{\hbox{Theorem \ref{L11-1}}}{\Longleftrightarrow}$ & $R$ is minimal \\
\hline
$\w_\pp=\w$ & $\stackrel{\hbox{Proposition \ref{coron-21}}}{\Longrightarrow}$ & $\la\pp\ra=\la R \ra$, $\pp$ is connected \\
\hline
$R$ finite & $\stackrel{\hbox{Popositon \ref{T3}}}{\Longrightarrow}$ &$\mm_c=\mm_m=\mm_r$, $c=$rank$(R)$\\
\hline
$\pp\in\mm_c$, $\mm_m$ or $\mm_r$ & $\stackrel{\hbox{Corollary \ref{T1}}}{\Longrightarrow}$& $\pp$ is finite  \\
\hline
$
X=A_1,\; B_\ell,\; C_\ell,\; F_4,\; G_2$
&$\stackrel{\hbox{Proposition \ref{T2}}}{\Longrightarrow}$&
$\begin{array}{l}
\hbox{all elements of}\;\mm_r\\
\hbox{have the same cardinality}
\end{array}$
\\
\hline
$X$= simply laced of rank$>1$&
$\underset{\text{Example \ref{ex2}}}{\overset{\text{Proposition \ref{min AD}}}{\Longrightarrow}}$
&$\mm_r=\mm_{rm}$ {and  $\mm_{rc}\subsetneq\mm_{r}$}\\
\hline
$\pp\in\mm_r$&$\stackrel{\hbox{Proposition \ref{min-gen.}}}{\Longrightarrow}$&$\begin{array}{l}
\pp\hbox{ is a minimal generating set for the}\\
 \hbox{conjugation presented group $\widehat\w$}\\
\end{array}$\\
\hline
$\mm_r=\mm_{rm}$&$\stackrel{\hbox{Theorem \ref{pro13}}}{\Longleftrightarrow}$& $R$ is minimal \\
\hline
$\begin{array}{l}
\nu\leq 2,\hbox{ or}\\
X=\hbox{simply laced of rank$>1$,}\\
X= F_4,\;G_2,\\
{X= A_1, B_\ell, C_\ell,\; \hbox{and}\; \ind(R)=0}\\
\end{array}$
&$\stackrel{\hbox{Corolary \ref{cor13}}}{\Longrightarrow}$& $\mm_r=\mm_{rm}$ \\
\hline
{$X=F_4, G_2$}&$\stackrel{\hbox{Corolary \ref{cor13}}}{\Longrightarrow}$&{$\mm_{r}=\mm_{rc}=\mm_{rm}$}\\
\hline
{$\begin{array}{c}
\ind(R)=0,\;X=A_1,\;B_\ell,\; C_\ell,\;\end{array}$}&
$\stackrel{\hbox{Corolary \ref{cor13}}}{\Longrightarrow}$&{$\mm_{r}=\mm_{rc}=\mm_{rm}$}\\
\hline
\end{tabular}}
\end{table}
}

\section{Relations to extended affine Lie algebras}\label{reala}
We begin with a brief introduction to the definition of an extended affine Lie algebra and some of its basic properties. 
For more details on the subject and further terminologies We refer the reader to \cite[Chapter I]{AABGP}. 

\subsection{Extended affine Lie algebras}\label{eeala}
\begin{DEF}\label{EALA}
	A triple ($\ll,\fm,\hh)$ in which $\ll$ is complex Lie algebra, $\hh$ is a subalgebra and $\fm$ is a bilinear form on $\ll$ is called an {\it extended affine Lie algebra} if the following 5 axioms hold:
	
	(E1) The form $\fm$ is symmetric, invariant, and non-degenerate.
	
	(E2) $\hh$ is non-trivial, finite-dimensional, and self-centralizing such that
	$$\ll=\sum_{\a\in R}\ll_\a$$
	where
	$$\ll_\a=\{x\in\ll\mid [h,x]=\a(h)(x)\hbox{ for all }h\in\hh\}$$
	and
	$$R=\{\a\in\hh^\star\mid\ll_\a\not=\{0\}\}.$$
	$R$ is called the {\it root system} of $\ll$ with respect to $\hh$. The form on $\ll$ restricted to $\hh$ is non-degenerate so 
	the form can be transferred to $\hh^\star$ in a natural way.
	We denote by $t_\a$ the unique element in $\hh$ which represents $\a\in\hh^\star$ via the form, namely
	$$\a(h)=(t_\a,h),\qquad(h\in\hh).
	$$
	Set
	$$R^0:=\{\a\in R\mid (\a,\a)=0\}\andd R^\times:=R\setminus R^0.$$
	
	(E4) $\ad(x)$ is locally nilpotent for $x\in\ll_\a$, $\a\in R^\times$.
	
	(E5) $R$ is irreducible, meaning that $R^\times$ is connected, and
	for $\sg\in R^0$, there exists $\a\in R^\times$ with $\a+\sg\in R$.	
\end{DEF}

\begin{pgraph}\label{b1}
	We assume from now on that $(\ll,\fm,\hh)$ is an extended affine Lie algebra with root system $R$.  One knows that $R$ is an extended affine root system in $\V:=\span_\bbbr R$, all non-isotropic root spaces are $1$-dimensional, and that for $\a\in R^\times$, one may choose $e_{\pm\a}\in\ll_{\pm\a}$ such that $(e_\a,h_\a:=[e_\a,e_{-\a}],e_{-\a})$ is an $\mathfrak{sl}_2$-triple.
	From now on, {\it we assume that $R$ is reduced}.  
	For $\a\in R^\times$ and $x_{\pm\a}\in\ll_{\pm\a}$, we define
	$\Phi_\a\in\Aut(\ll)$ by
	$$\Phi_\a:=\exp(\ad x_\a)\exp(-\ad x_{-\a})\exp(\ad x_\a).$$
	Then $\Phi_\a(\ll_\be)=\ll_{w_\a(\be)}$ for $\be\in R$, see \cite[Proposition I.1.27]{AABGP}.
	{As in Definition \ref{eawg}, we denote by $\w$ the extended affine Weyl group of $R$. We also denote by $\w_\ll$ the Weyl group of $\ll$, the subgroup of $\GL(\hh^\star)$ generated by reflections $\be\mapsto \be-(\be,\a^\vee)\a$, $\a\in R^\times$.}
	Then $\Phi_\a$ restricted to $\hh\equiv\hh^\star$ coincides with the reflection $w_\a\in\w$, identifying $\w_\ll$ with the extended affine Weyl group $\w$ of $R$. 
	
	We recall that the core $\ll_c$ of $\ll$ is the subalgebra of $\ll$ generated by the non-isotropic root spaces $\ll_\a$, $\a\in R^\times$.
\end{pgraph}

\subsection{Reflectable bases and the core}
We begin with a lemma that shows, using the concept of a reflectable set, that the root spaces associated with a reflectable set (up to plus-minus sign) generate the core, implying that the core of an extended affine Lie algebra is finitely generated (compare with \cite[Section 6.12]{Neh} and \cite[Corollary 2.14]{A4}). We start with terminology.

\begin{DEF}\label{defl1}
	We call a subset $\pp$ of $R^\times$, a {\it root-generating set} for $\ll$ if the root spaces $\ll_{\pm\a}$, $\a\in\pp$ generate $\ll_c$.
	We call a {root-generating set} {\it minimal} if no proper subset of $\pp$ is a root-generating set.
\end{DEF}

\begin{lem}\label{llem1}
	Any reflectable set for $R$ is a root-generating set for $\ll$. In particular, $\ll_c$ is finitely generated.
\end{lem}

\proof
Let $\pp$ be a reflectable set and $\a\in R^\times.$  Then $\a=w_{\a_1}\cdots w_{\a_t}(\a_{t+1})$, for some $\a_i\in\pp$. Therefore,
$$\Phi_{\a_{1}}\cdots\Phi_{\a_t}(\ll_{\a_{t+1}})=
\ll_{w_{\a_1}\cdots w_{\a_t}(\a_{t+1})}=\ll_\a.$$
{Thus $\ll_\a$ contains in the subalgebra of $\ll$ generated by $\ll_\be$, $\be\in\pp^\pm$.} The second statement now follows from Corollary \ref{T1}.\qed

\begin{lem}\label{llem2}
	Let $R$ be a simply laced extended affine root system of type $X$ and rank $>1$. Let $\pp$ 
	be a reflectable base for $R$. Then $\pp$ is a minimal root-generating set for $\ll$. In particular the set $\pp(X)$ given in Table \ref{tab:fq} is a minimal root-generating set for
	$\ll$.
\end{lem}

\proof
Assume to the contrary that $\pp$ contains a proper root-generating set $\pp'$.
If $\a\in\pp$, then $\ll_\a$ is generated by some root spaces
$\ll_{\a_1},\ldots,\ll_{\a_t}$ for some $\a_1,\ldots,\a_t\in\pp'$.
Thus $\a\in\la\pp'\ra$, and so $\la\pp'\ra=\la\pp\ra=\la R\ra.$
Then by Theorem {3.2} of \cite{ASTY1}, $\pp'$ is a reflectable set for $R$. But this contradicts the minimality of $\pp$ as a reflectable base.
The second assertion in the statement follows from Proposition \ref{L11}.\qed

\begin{lem}\label{llem3}
	If $\ind(R)=0$, then any reflectable base $\pp$ is a minimal root-generating set for $\ll$. In particular, if $R$ is simply laced of rank $>1$, or is of type $F_4$ or $G_2$, then any reflectable base is a minimal root-generating set for $\ll$.\qed
\end{lem}

\proof If $R$ is simply laced of rank $>1$, then we are done by Lemma {\ref{llem2}}. For the remaining types, we have from Proposition \ref{T2} that all reflectable bases have the same cardinality. Now if $\ind(R)=0$ then we see from {\pref{note2} that $|\pp|=|\pp(X)|=\ell+\nu$  which is equal to the rank} of the free abelian group $\la R\ra$. Since the roots corresponding to any generating set of root vectors for $\ll_c$ must generate $\la R\ra$, the first assertion in the statement holds. The second assertion is now clear as $\ind(R)=0$ for the given types.\qed 

The following example shows that if $\ind(R)\not=0$ then Lemma \ref{llem3} may fail. Also, it shows that an extended affine Lie algebra might have minimal root-generating sets with different cardinalities.

\begin{exa}\label{ex23}
	Suppose $\ll$ is an extended affine Lie algebra with root system $R$. 
	
	(i) Let $R$ be of type $A_1$ and nullity $2$ with  $\ind(R)>0$. We know that (see \pref{note2}) up to isomorphism
	$R=\Lam+(\pm\a+\Lam)$ where $\Lam=\bbbz\sg_1\oplus\bbbz\sg_2$.
	We consider the reflectable base
	$\pp(X)=\{\a,\sg_1-\a,\sg_2-\a,\sg_1+\sg_2-\a\}$ for $R$, see Table \ref{tab:fq}.
	From \cite[Remark 1.5(ii)]{A5}, we know that
	$$[[\ll_{\sg_1-\a},\ll_\a],\ll_{\sg_2-\a}]]\not=\{0\}.$$
	Thus the $1$-dimensional space $\ll_{\sg_1+\sg_2-\a}$ is generated by root spaces corresponding to roots $\{\a,\sg_1-\a,\sg_2-\a\}$.
	Therefore by Lemma \ref{llem1}, $\ll_c$ is generated by root spaces
	$\ll_{\pm\be}$, $\be\in\{\a,\sg_1-\a,\sg_2-\a\}$. Thus $\pp(X)$ is not a minimal root-generating set for $\ll_c$.
	
	(ii) Let $R$ be of type $A_2$. Then up to isomorphism
	$R=\Lam\cup (\dot R+\Lam)$ where $\Lam=\bbbz\sg_1\oplus\bbbz\sg_2$ and $\dot R$ is a finite root system of type $A_2$ with a base $\{\a_1,\a_2\}$.
	Then $\pp=\{\a_1,\a_2,\a_1+\sg_1,\a_2+\sg_2\}$ and $\pp'=\{\a_1,\a_2,\a_1+2\sg_1,\a_2+3\sg_2,\a_2+\sg_2\}$ are reflectable bases for $R$, see Example \ref{ex2}. By Lemma \ref{llem3}  both $\pp$ and $\pp'$ are minimal root-generating sets with $|\pp|\not=|\pp'|$ for the corresponding Lie algebra $\ll$.
\end{exa}

\subsection{Elliptic Lie algebras versus $2$-extended affine Lie algebras}
We begin with a definition.
By a {\it $2$-extended affine Lie algebra (root system)} we mean an extended affine Lie algebra (root system) of nullity $2$. Throughout this section, $R$ is a {reduced} $2$-extended affine root system.

\begin{DEF}\label{elliptic}
	(i) In the literature, a $2$-extended affine root system is referred to as an {\it elliptic root system}.
	
	(ii) 
	Let $\ll$ be a Lie algebra containing a subalgebra $\hh$ with a decomposition $\ll=\sum_{\a\in\hh^\star}\ll_\a$, where
	$$\ll_\a=\{x\in\ll\mid [h,x]=\a(h)x\hbox{ for all }h\in\hh\}.
	$$
	The set $R:=\{\a\in\hh^\star\mid\ll_\a\not=\{0\}\}$ is called
	the {\it root system} of $\ll$ with respect to $\hh$. We call $(\ll,\hh)$, or simply $\ll$  an {\it elliptic Lie algebra} if $R$ is an elliptic root system, with respect to a symmetric positive semi-definite form on the real span of $R$.
\end{DEF}

\begin{cor}\label{corl1}
	Let $\ll$ be a $2$-extended affine  Lie algebra of type $X$, rank $\ell$ and twist number $t$. Let $S$ be as in \pref{pgraph2}. Suppose either of the following holds:
	
	- $X$ is simply laced of rank$>1$,
	
	- $X=A_1$ and $\ind(S)=2$,
	
	- $X=F_4,\;G_2$,
	
	- $X=B_2$ with $t=1$; or  $t=0,2$ and $\ind(S)=2$,
	
	- $X=B_\ell$ ($\ell>2$) with $t=0,1$; or $t=2$ and $\ind(S)=2$,
	
	- $X=C_\ell$ with $t=1,2$; or $t=0$ and $\ind(S)=2$.
	
	Then any reflectable base is a minimal root-generating set for  $\ll$.
\end{cor}

\proof For all the types given in the statement we have $\ind(R)=0$ and so by Lemma \ref{llem3} we are done.\qed

\subsection{A prototype presentation}
{Let 
	$R=R^0\cup(\dot{R}_{sh}+S)\cup (\rdl+L)$ be an elliptic root system,
	where $\rd=\rds\cup\rdl\cup\{0\}$ and $S,L$ are defined in \pref{pgraph2}.} Let
$\dot{\pp}=\{\a_1, \ldots, \a_\ell\}$ be a fundamental system for $\rd$ and {consider the Dynkin diagrams as follow. }  
\hspace{.5cm}
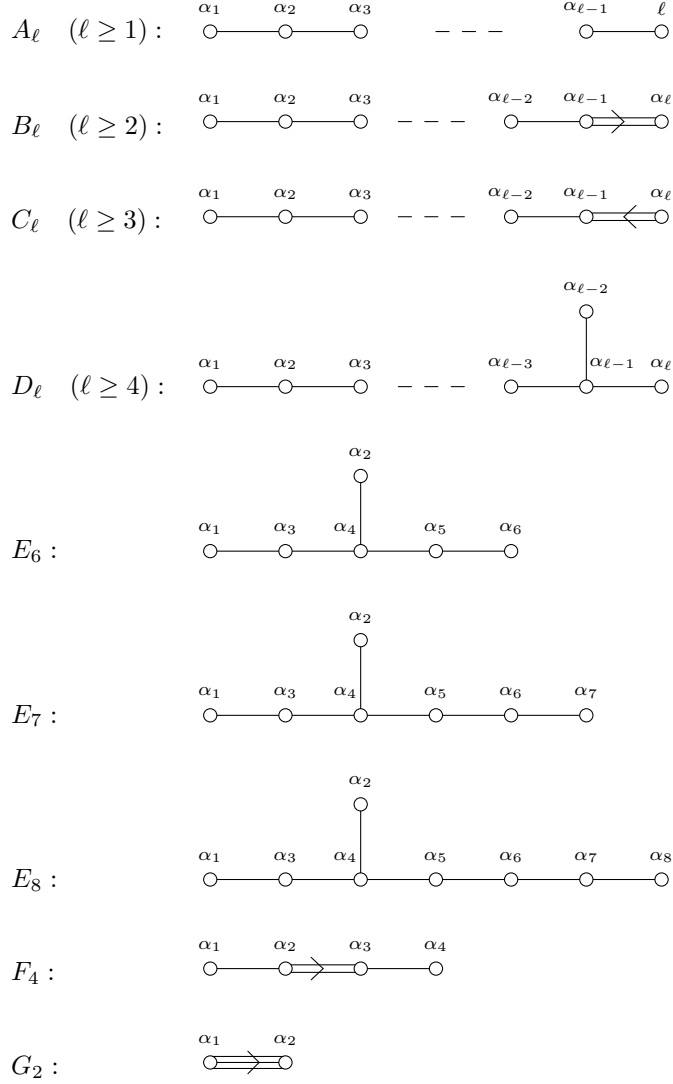
$$\begin{array}{ll}
	\vspace{.5cm}
	A_\ell\quad(\ell\geq 1):&
	\definecolor{qqqqff}{rgb}{0.,0.,1.}
	\definecolor{siah}{rgb}{1.,1.,1.}
	\begin{tikzpicture}[line cap=round,line join=round,>=triangle 45,x=1.0cm,y=1.0cm]
		\draw (-6.,3.)-- (-5.,3.);
		\draw (-5.,3.)--(-4.,3.);
		\draw [dash pattern=on 5pt off 5pt] (-3.,3.)-- (-2.,3.);
		\draw (-1.,3.)-- (0.,3.);
		\begin{scriptsize}
			\draw [fill=white] (-6.,3.) circle (2.5pt);
			\draw[color=black] (-6.,3.3) node {$\a_1$};
			\draw [fill=white] (-5.,3.) circle (2.5pt);
			\draw[color=black] (-5.,3.3) node {$\a_2$};
			\draw [fill=white] (-4.,3.) circle (2.5pt);
			\draw[color=black] (-4.,3.3) node {$\a_3$};
			\draw [fill=white] (-1.,3.) circle (2.5pt);
			\draw[color=black] (-1.,3.3) node {$\a_{\ell-1}$};
			\draw[color=black] (0.,3.3) node {$\ell$};
			\draw [fill=white] (0.,3.) circle (2.5pt);
		\end{scriptsize}
	\end{tikzpicture}\\
	\vspace{.5cm}	B_\ell\quad(\ell\geq 2):&
	\definecolor{qqqqff}{rgb}{0.,0.,1.}
	\definecolor{siah}{rgb}{1.,1.,1.}
	\begin{tikzpicture}[line cap=round,line join=round,>=triangle 45,x=1.0cm,y=1.0cm]
		\draw (-6.,3.)-- (-5.,3.);
		\draw (-5.,3.)--(-4.,3.);
		\draw (-2.,3.)--(-1.,3.);
		\draw [dash pattern=on 5pt off 5pt] (-3.5,3.)-- (-2.5,3.);
		\draw (-1.,3.05)-- (0.,3.05);
		\draw (-1.,2.95)-- (0.,2.95);
		\draw (-.5,3.)-- (-.66,3.15);
		\draw (-.5,3.)-- (-.66,2.85);
		\begin{scriptsize}
			\draw [fill=white] (-6.,3.) circle (2.5pt);
			\draw[color=black] (-6.,3.3) node {$\a_1$};
			\draw [fill=white] (-5.,3.) circle (2.5pt);
			\draw[color=black] (-5.,3.3) node {$\a_2$};
			\draw [fill=white] (-4.,3.) circle (2.5pt);
			\draw[color=black] (-4.,3.3) node {$\a_3$};
			\draw [fill=white] (-2.,3.) circle (2.5pt);
			\draw[color=black] (-2.,3.3) node {$\a_{\ell-2}$};
			\draw [fill=white] (-1.,3.) circle (2.5pt);
			\draw[color=black] (-1.,3.3) node {$\a_{\ell-1}$};
			\draw[color=black] (0.,3.3) node {$\a_{\ell}$};
			\draw [fill=white] (0.,3.) circle (2.5pt);
		\end{scriptsize}
	\end{tikzpicture}\\
	\vspace{.5cm}	C_\ell\quad(\ell\geq 3):&
	\definecolor{qqqqff}{rgb}{0.,0.,1.}
	\definecolor{siah}{rgb}{1.,1.,1.}
	\begin{tikzpicture}[line cap=round,line join=round,>=triangle 45,x=1.0cm,y=1.0cm]
		\draw (-6.,3.)-- (-5.,3.);
		\draw (-5.,3.)--(-4.,3.);
		\draw (-2.,3.)--(-1.,3.);
		\draw [dash pattern=on 5pt off 5pt] (-3.5,3.)-- (-2.5,3.);
		\draw (-1.,3.05)-- (0.,3.05);
		\draw (-1.,2.95)-- (0.,2.95);
		\draw (-.5,3.)-- (-.34,3.15);
		\draw (-.5,3.)-- (-.34,2.85);
		\begin{scriptsize}
			\draw [fill=white] (-6.,3.) circle (2.5pt);
			\draw[color=black] (-6.,3.3) node {$\a_1$};
			\draw [fill=white] (-5.,3.) circle (2.5pt);
			\draw[color=black] (-5.,3.3) node {$\a_2$};
			\draw [fill=white] (-4.,3.) circle (2.5pt);
			\draw[color=black] (-4.,3.3) node {$\a_3$};
			\draw [fill=white] (-2.,3.) circle (2.5pt);
			\draw[color=black] (-2.,3.3) node {$\a_{\ell-2}$};
			\draw [fill=white] (-1.,3.) circle (2.5pt);
			\draw[color=black] (-1.,3.3) node {$\a_{\ell-1}$};
			\draw[color=black] (0.,3.3) node {$\a_{\ell}$};
			\draw [fill=white] (0.,3.) circle (2.5pt);
		\end{scriptsize}
	\end{tikzpicture}\\
	\vspace{.5cm}D_\ell\quad(\ell\geq 4):&
	\definecolor{qqqqff}{rgb}{0.,0.,1.}
	\definecolor{siah}{rgb}{1.,1.,1.}
	\begin{tikzpicture}[line cap=round,line join=round,>=triangle 45,x=1.0cm,y=1.0cm]
		\draw (-6.,3.)-- (-5.,3.);
		\draw (-5.,3.)--(-4.,3.);
		\draw (-2.,3.)--(-1.,3.);
		\draw [dash pattern=on 5pt off 5pt] (-3.5,3.)-- (-2.5,3.);
		\draw (-1.,3.)-- (-1.,4.);
		\draw (-1.,3.)-- (0.,3.);
		\begin{scriptsize}
			\draw [fill=white] (-6.,3.) circle (2.5pt);
			\draw[color=black] (-6.,3.3) node {$\a_1$};
			\draw [fill=white] (-5.,3.) circle (2.5pt);
			\draw[color=black] (-5.,3.3) node {$\a_2$};
			\draw [fill=white] (-4.,3.) circle (2.5pt);
			\draw[color=black] (-4.,3.3) node {$\a_3$};
			\draw [fill=white] (-2.,3.) circle (2.5pt);
			\draw[color=black] (-2.,3.3) node {$\a_{\ell-3}$};
			\draw [fill=white] (-1.,3.) circle (2.5pt);
			\draw[color=black] (-1.,4.3) node {$\a_{\ell-2}$};
			\draw [fill=white] (-1.,4.) circle (2.5pt);
			\draw[color=black] (-0.65,3.3) node {$\a_{\ell-1}$};
			\draw [fill=white] (0.,3.) circle (2.5pt);
			\draw[color=black] (0.,3.3) node {$\a_{\ell}$};
		\end{scriptsize}
	\end{tikzpicture}\\
	\vspace{.5cm}E_6:&
	\definecolor{qqqqff}{rgb}{0.,0.,1.}
	\definecolor{siah}{rgb}{1.,1.,1.}
	\begin{tikzpicture}[line cap=round,line join=round,>=triangle 45,x=1.0cm,y=1.0cm]
		\draw (-6.,3.)-- (-5.,3.);
		\draw (-5.,3.)--(-4.,3.);
		\draw (-4.,3.)--(-3.,3.);
		\draw (-3.,3.)-- (-2.,3.);
		\draw (-4.,3.)-- (-4.,4.);
		\begin{scriptsize}
			\draw [fill=white] (-6.,3.) circle (2.5pt);
			\draw[color=black] (-6.,3.3) node {$\a_1$};
			\draw [fill=white] (-5.,3.) circle (2.5pt);
			\draw[color=black] (-5.,3.3) node {$\a_3$};
			\draw [fill=white] (-4.,3.) circle (2.5pt);
			\draw[color=black] (-4.2,3.3) node {$\a_4$};
			\draw [fill=white] (-3.,3.) circle (2.5pt);
			\draw[color=black] (-3.,3.3) node {$\a_5$};
			\draw [fill=white] (-2.,3.) circle (2.5pt);
			\draw[color=black] (-2.,3.3) node {$\a_6$};
			\draw [fill=white] (-4.,4.) circle (2.5pt);
			\draw[color=black] (-4.,4.3) node {$\a_2$};
		\end{scriptsize}
	\end{tikzpicture}\\
	\vspace{.5cm}E_7:&
	\definecolor{qqqqff}{rgb}{0.,0.,1.}
	\definecolor{siah}{rgb}{1.,1.,1.}
	\begin{tikzpicture}[line cap=round,line join=round,>=triangle 45,x=1.0cm,y=1.0cm]
		\draw (-6.,3.)-- (-5.,3.);
		\draw (-5.,3.)--(-4.,3.);
		\draw (-4.,3.)--(-3.,3.);
		\draw (-2.,3.)-- (-1.,3.);
		\draw (-3.,3.)-- (-2.,3.);
		\draw (-4.,3.)-- (-4.,4.);
		\begin{scriptsize}
			\draw [fill=white] (-6.,3.) circle (2.5pt);
			\draw[color=black] (-6.,3.3) node {$\a_1$};
			\draw [fill=white] (-5.,3.) circle (2.5pt);
			\draw[color=black] (-5.,3.3) node {$\a_3$};
			\draw [fill=white] (-4.,3.) circle (2.5pt);
			\draw[color=black] (-4.2,3.3) node {$\a_4$};
			\draw [fill=white] (-3.,3.) circle (2.5pt);
			\draw[color=black] (-3.,3.3) node {$\a_5$};
			\draw [fill=white] (-2.,3.) circle (2.5pt);
			\draw[color=black] (-2.,3.3) node {$\a_6$};
			\draw[color=black] (-1.,3.3) node {$\a_7$};
			\draw [fill=white] (-1.,3.) circle (2.5pt);
			\draw [fill=white] (-4.,4.) circle (2.5pt);
			\draw[color=black] (-4.,4.3) node {$\a_2$};
		\end{scriptsize}
	\end{tikzpicture}\\
				\vspace{.5cm}E_8:&
				\definecolor{qqqqff}{rgb}{0.,0.,1.}
				\definecolor{siah}{rgb}{1.,1.,1.}
				\begin{tikzpicture}[line cap=round,line join=round,>=triangle 45,x=1.0cm,y=1.0cm]
					\draw (-6.,3.)-- (-5.,3.);
					\draw (-5.,3.)--(-4.,3.);
					\draw (-4.,3.)--(-3.,3.);
					\draw (-2.,3.)-- (-1.,3.);
					\draw (-1.,3.)--(0.,3.0);
					\draw (-3.,3.)-- (-2.,3.);
					\draw (-4.,3.)-- (-4.,4.);
					\begin{scriptsize}
						\draw [fill=white] (-6.,3.) circle (2.5pt);
						\draw[color=black] (-6.,3.3) node {$\a_1$};
						\draw [fill=white] (-5.,3.) circle (2.5pt);
						\draw[color=black] (-5.,3.3) node {$\a_3$};
						\draw [fill=white] (-4.,3.) circle (2.5pt);
						\draw[color=black] (-4.2,3.3) node {$\a_4$};
						\draw [fill=white] (-3.,3.) circle (2.5pt);
						\draw[color=black] (-3.,3.3) node {$\a_5$};
						\draw [fill=white] (-2.,3.) circle (2.5pt);
						\draw[color=black] (-2.,3.3) node {$\a_6$};
						\draw[color=black] (-1.,3.3) node {$\a_7$};
						\draw [fill=white] (-1.,3.) circle (2.5pt);
						\draw [fill=white] (-4.,4.) circle (2.5pt);
						\draw[color=black] (-4.,4.3) node {$\a_2$};
						\draw [fill=white] (0.,3.) circle (2.5pt);
						\draw[color=black] (0.,3.3) node {$\a_8$};
					\end{scriptsize}
				\end{tikzpicture}\\
				
					\vspace{.5cm}F_4:&
					\definecolor{qqqqff}{rgb}{0.,0.,1.}
					\definecolor{siah}{rgb}{1.,1.,1.}
					\begin{tikzpicture}[line cap=round,line join=round,>=triangle 45,x=1.0cm,y=1.0cm]
						\draw (-5.,3.)--(-4.,3.);
						\draw (-2.,3.)--(-3.,3.);
						\draw (-3.,3.05)-- (-4.,3.05);
						\draw (-3.,2.95)-- (-4.,2.95);
						\draw (-3.5,3.)-- (-3.66,3.15);
						\draw (-3.5,3.)-- (-3.66,2.85);
						\begin{scriptsize}
							\draw [fill=white] (-5.,3.) circle (2.5pt);
							\draw[color=black] (-5.,3.3) node {$\a_1$};
							\draw [fill=white] (-4.,3.) circle (2.5pt);
							\draw[color=black] (-4.,3.3) node {$\a_2$};
							\draw [fill=white] (-3.,3.) circle (2.5pt);
							\draw[color=black] (-3.,3.3) node {$\a_3$};
							\draw [fill=white] (-2.,3.) circle (2.5pt);
							\draw[color=black] (-2.,3.3) node {$\a_4$};
						\end{scriptsize}
					\end{tikzpicture}\\
					\vspace{.5cm}G_2:&
					\definecolor{qqqqff}{rgb}{0.,0.,1.}
					\definecolor{siah}{rgb}{1.,1.,1.}
					\begin{tikzpicture}[line cap=round,line join=round,>=triangle 45,x=1.0cm,y=1.0cm]
						\draw (-6.,3.08)-- (-5.,3.08);
						\draw (-6.,3.)-- (-5.,3.);
						\draw (-6.,2.92)-- (-5.,2.92);
						\draw (-5.35,3.)-- (-5.5,3.15);
						\draw (-5.35,3.)-- (-5.5,2.85);
						\begin{scriptsize}
							\draw [fill=white] (-5.,3.) circle (2.5pt);
							\draw[color=black] (-5.,3.3) node {$\a_2$};
							\draw [fill=white] (-6.,3.) circle (2.5pt);
							\draw[color=black] (-6.,3.3) node {$\a_1$};
						\end{scriptsize}
					\end{tikzpicture}\\
				\end{array}$$

				As in \pref{e1}, we set
				$\V=\span_\bbbr R$, $
				\wtilde{\V}=\V\oplus(\V^0)^\star$ and we let  $\fm$ be the non-degenerate form on $\wtilde\V$. We recall that $\V=\dot\V\oplus\V^0$, and that  $\V^0=\bbbr\sg_1\oplus\bbbr\sg_2$, where $\sg_1,\sg_2\in R^0$ and $\la R^0\ra=\bbbz\sg_1\oplus\bbbz\sg_2$. We fix a basis $\{\lam_1,\lam_2\}$ for $(\V^0)^\star$ by
				\begin{equation}
					\label{fix1}
					\lam_i(\sg_j)=\d_{ij},\quad(1\leq i,j\leq 2).
				\end{equation}
				We consider $\lam_1,\lam_2$ as elements of $\V^\star$ by requiring $\lam_1(\dot\V)=
				\lam_2(\dot\V)=\{0\}$. 
				
				{Assume that $\pp=\pp(X)$ is} the reflectable base for $R$ given in Table \ref{ttt1}. 	
				{\scriptsize
					\begin{table}[ht]
						\caption{A reflectable base {for the elliptic root system $R$}}\label{ttt1}
						\vspace{-4mm}
						\begin{tabular}
							[c]{|c | c |} \hline
							
							$\mbox{Type}$& $\pp(X)$\\
							\whline
							$A _1$ & $\begin{array}{c}
								\{\a_1, \sg_1 -\a_1,\sg_2-\a_1\}\hbox{ if }\ind(R)=0\\
								\{\a_1, \sg_1 -\a_1,\sg_2-\a_1,\sg_1+
								\sg_2-\a_1\}\hbox{ if }\ind(R)=1
							\end{array}$\\
							\hline
							$\begin{array}{c}
								A_\ell,\; (\ell>1) \\
								D_\ell\\
								E_{6.7.8}
							\end{array}$ & $\{\a_1, \ldots, \a_\ell, \sg_1-\a_1, \sg_2-\a_1\}$\\
							\hline
							$F_4$ & $\begin{array}{c}\{\a_1, \a_2,
								\a_3,\a_4,  \sg_{1}-\a_1,\sg_2 -\a_1\}\hbox{ if }t=0\\
								\{\a_1, \a_2,
								\a_3,\a_4,  \sg_{1}-\a_4,\sg_2 -\a_1\}\hbox{ if }t=1\\
								\{\a_1, \a_2,
								\a_3,\a_4,  \sg_{1}-\a_4,\sg_2 -\a_4\}\hbox{ if }t=2
							\end{array}
							$\\
							\hline
							$G_2$ & $\begin{array}{c}\{\a_1, \a_2,
								\sg_{1}-\a_1,\sg_2 -\a_1\}\hbox{ if }t=0\\
								\{\a_1, \a_2,\sg_{1}-\a_2,\sg_2 -\a_1\}\hbox{ if }t=1\\
								\{\a_1, \a_2,
								\sg_{1}-\a_2,\sg_2 -\a_2\}\hbox{ if }t=2
							\end{array}
							$\\
							\hline
							$B_2$ & $\begin{array}{c}\{\a_1, \a_2,  \sg_{1}-\a_1, \sg_2 -\a_1 \}\hbox{ if }t=0,\;\ind(R)=0\\
								\{\a_1, \a_2,  \sg_{1}-\a_1, \sg_2 -\a_1,\sg_1+\sg_2-\a_1 \}\hbox{ if }t=0,\;\ind(R)=1\\	
								\{\a_1, \a_2,  \sg_{1}-\a_1, \sg_2 -\a_2 \}\hbox{ if }t=1
								\\
								\{\a_1, \a_2,  \sg_{1}-\a_2, \sg_2 -\a_2 \}\hbox{ if }t=2,\;\ind(R)=0,\\
								\{\a_1,\a_2,  \sg_{1}-\a_2, \sg_2 -\a_2,\sg_1+\sg_2-\a_{2} \}\hbox{ if }t=2,\;\ind(R)=1
							\end{array}
							$\\
							\hline
							$B_\ell,\; (\ell>2)$ & $\begin{array}{c}\{\a_1, \ldots, \a_\ell,  \sg_{1}-\a_1, \sg_2 -\a_1 \}\hbox{ if }t=0\\
								\{\a_1, \ldots, \a_\ell,  \sg_{1}-\a_1, \sg_2 -\a_\ell \}\hbox{ if }t=1
								\\
								\{\a_1, \ldots, \a_\ell,  \sg_{1}-\a_\ell, \sg_2 -\a_\ell \}\hbox{ if }t=2,\;\ind(R)=0,\\
								\{\a_1, \ldots, \a_\ell,  \sg_{1}-\a_\ell, \sg_2 -\a_\ell,\sg_1+\sg_2-\a_{\ell} \}\hbox{ if }t=2,\;\ind(R)=1
							\end{array}
							$\\
							\hline
							$C_\ell,\; (\ell>2)$ & $\begin{array}{c}\{\a_1, \ldots, \a_\ell,  \sg_{1}-\a_1, \sg_2 -\a_1 \}\hbox{ if }t=2\\
								\{\a_1, \ldots, \a_\ell,  \sg_{1}-\a_1, \sg_2 -\a_\ell \}\hbox{ if }t=1
								\\
								\{\a_1, \ldots, \a_\ell,  \sg_{1}-\a_\ell, \sg_2 -\a_\ell \}\hbox{ if }t=0,\;\ind(R)=0,\\
								\{\a_1, \ldots, \a_\ell,  \sg_{1}-\a_\ell, \sg_2 -\a_\ell,\sg_1+\sg_2-\a_{\ell} \}\hbox{ if }t=0,\;\ind(R)=1
							\end{array}
							$\\
							\hline
						\end{tabular}
					\end{table}   
				}
			{We recall from Proposition \ref{coron2} and axiom (R4) of Definition \ref{def2} that
				$$R=
				\w_\pp\pp\cup\big((\w_\pp\pp-\w_\pp\pp)\cap\V^0\big),
				$$
				so $R$ is uniquely determined by $\pp$.}

			For $\a,\be\in R^\times$, we define $$n_{\a,\be}:=\hbox{min}\{n\in\bbbz_{>0}\mid n\a+\be\not\in R\}.$$
			
			In what follows by $\pp^\pm$ we mean $\pp\cup(-\pp)$. Also a term
			$[x_1,\ldots, x_n]$ in a Lie algebra means
			$[x_1,[x_2,[x_3,\ldots,[x_{n-1},x_n]\cdots]]].$
			
			Let $\what\ll$ be the Lie algebra defined by generators
			$$\what{X}_{\a},\;\what{H}_\a,
			\what{d}_1,\;\what d_2,\qquad(\a\in\pp^\pm),
			$$
			subject to the following relations
				\begin{equation}\label{eq111}
					\begin{array}{ll}
						\hbox{(I)}&\sum_{i=1}^k\what{H}_{\be_i}=\what{H}_{\sum_{i=1}^k\be_i}\hbox { if }\be_1,\ldots,\be_{k},\sum_{i=1}^k\be_i\in\pp,\\
						
						\hbox{(II)}& [\what{H}_\a,\what{H}_\be]=0,\; \a,\be\in\pp^\pm,\\
						
						
						\hbox{(III)}& [\what{H}_\a,\what{X}_\be]=(\be,\a^\vee)\what{X}_\be,\;\a,\be\in\pp^\pm,\\
						
						\hbox{(IV)}&[\what{X}_\a,\what{X}_{-\a}]=\what{H}_\a,\;\a\in\pp^\pm,\\
						
						{\hbox{(V)}}&{[\what{X}_{\a_1},\ldots,\what X_{\a_n}]=0,\;\a_i\in\pp^\pm,\sum_{i=1}^n\a_i\not\in R,}
						\\
						
						
						\hbox{(VI)}&[\what d_i,\what X_\a]=\lam_i(\a)\what X_\a,\quad 
						[\what d_i,\what d_j]=[\what d_i,\what H_\a]=0,\quad \a\in\pp^\pm,\; 1\leq i,j\leq 2.
					\end{array}
				\end{equation}
			We see from defining relations of the form (V) that
			
			{$\hbox{(V)}'\quad\ad\what{X}_{\a}^{n_{\a,\be}}(\what{X}_{\be})=0, \;\a,\be\in\pp^\pm.$}
			
			{
				\begin{rem}\label{revrem2}
					Let $\Phi$ be an irreducible finite or affine root system, and $\pp$ be a fundamental system for $\Phi$. Let $\what\ll$ be the Lie algebra defined by  generators $\what{X}_\a, \what H_\a, \what d_1$, $\a\in{\pp}^\pm$, and relations (I)-(IV), (V$'$), (VI), (while $\what{d}_2$ is dropped). Then $\what\ll$ is either a finite-dimensional simple Lie algebra or the derived subalgebra of an affine Lie algebra with root system $\Phi$, see \cite{Kac}. Since relations of the form (V) hold in $\what\ll$, it follows that (V) and (V$'$) are equivalent in this case.
				\end{rem}
			}
			
			To proceed, we recall that in \cite[Chapter III]{AABGP}, associated to each elliptic root system $R$ an extended affine Lie algebra $(\ll,\fm,\hh)$ with root system $R$ is constructed such that
			$\ll=\ll_c\oplus{\mathbb C}d_1\oplus{\mathbb C}d_2$, where $d_1,d_2$ act as derivations on $\ll_c$, and $[d_1,d_2]=0$. {Moreover, 
				$\hh=\sum_{\dot\a\in\dot R}\bbbc t_{\dot\a}\oplus\bbbc t_{\sg_1}\oplus\bbbc t_{\sg_2}$, and $(d_i,t_{\sg_j})=\d_{ij}=\lam_i(\sg_j)$, $i,j=1,2$. So we may identify $\lam_i\in{(\V^0)}^\star\sub\hh^{\star\star}\equiv\hh$ with $d_i$. Then for $\a\in R$, $\lam_i(\a)=\a(d_i).$
				We fix this Lie algebra $\ll$ throughout the section and set
				$\hh_c:=\hh\cap\ll_c$.}
			
			For each $\a\in R^\times$, we fix an $\mathfrak{sl}_2$-triple
			$(e_\a,h_\a,e_{-\a})$, where $e_{\pm\a}\in\ll_{\pm\a}$,
			$[e_\a,e_{-\a}]=h_\a$ and $[h_\a,e_{\pm\a}]=\pm 2h_{\a}$, see \pref{b1}.

			Since relations (\ref{eq111}) hold in $\ll$ for $e_{\a}, h_\a, d_1,d_2$, $\a\in\pp^\pm$, in place of
			$\what X_\a$, $\what H_\a,\what d_1,\what d_2$, and since by Lemma \ref{llem1}
			the elements $e_{\a}, h_\a, d_1, d_2$, $\a\in\pp^\pm$ generate $\ll$, we get an epimorphism 
			\begin{equation}\label{eq1234}
				\begin{array}{c}
					\Psi:\what\ll\rightarrow\ll,\\
					\what{X}_{\a}\mapsto e_{\a},\quad \what{H}_\a\mapsto h_\a,\;(\a\in\pp^\pm),\quad\what d_i\mapsto d_i.\\
				\end{array}
			\end{equation}
			We set
			$$\what\hh:=\span_{\mathbb C}\{\what H_\a,\what d_1,\what d_2\mid\a\in\pp^\pm\}=
			\span_{\mathbb C}\{\what H_\a,\what d_1,\what d_2\mid\a\in\pp\}.$$
			We note that the second equality holds here, since $\what H_{-\a}=-\what H_\a$, for $\a\in\pp$, by relations of the form (IV) in (\ref{eq111}).
			
			\begin{lem}\label{llem16} $\Psi$ restricts to an isomorphism $\what\hh\rightarrow\hh$. In particular, $\what\hh$ is an abelian subalgebra of $\what\ll$ with
				$\dim\what\hh=\dim\hh=\ell+4$.
			\end{lem} 
			
			\proof Consider the epimorphism $\Psi:\what\ll\rightarrow \ll$ given by (\ref{eq1234}). From \cite[Chapter III]{AABGP}, we know that $\dim\hh=\ell+4$ and $\hh=\sum_{\a\in\pp}{\mathbb C} h_\a\oplus{\mathbb C} d_1\oplus{\mathbb C} d_2$. So if $|\pp=\pp(X)|=\ell+2$, then $\{h_\a, d_1, d_2\mid\a\in\pp\}$ is a linearly independent set and so its preimage $\{\what{H}_\a\mid\a\in\pp\}$ under $\Psi$ is also a linearly independent set; therefore  we are done in this case. 
			
			Assume next that  $|\pp|>\ell+2$.
			By Table \ref{ttt1}, we are in one of the cases 
			$X=A_1,$ $X=B_\ell$, $C_\ell$, with $\ind (R)=1$.
			For type $A_1$, we have $\pp(A_1)=\{\a_1,\sg_1-\a_1,\sg_2-\a_1,
			\sg_1+\sg_2-\a_1\}$. Now the elements $\what{H}_{\a_1}$, $\what{H}_{\sg_1-\a_1}$ and $\what{H}_{\sg_2-\a_1}$ are linearly independent in $\what\hh$ as their images are so in $\hh$. Also
			from relations of the form (I) in (\ref{eq111}), we have
			$\what{H}_{\sg_1+\sg_2-\a}=\what{H}_{\sg_1-\a}+\what{H}_{\sg_2-\a}+\what{H}_\a.$ Thus $\dim\what{\hh}=\dim\hh.$
			
			Assume next that $X=B_\ell$ $(\ell>2)$. From Table \ref{ttt1}, we see that the only  reflectable base $\pp$ with $|\pp|>\ell+2$ is the one given in the forth row, namely
			$$\{\a_1, \ldots, \a_\ell,  \sg_{1}-\a_\ell, \sg_2 -\a_\ell,\sg_1+\sg_2-\a_{\ell} \}.$$
			We know that $\what{H}_{\a_1},\ldots\what{H}_{\a_\ell},\what{H}_{\sg_1-\a_\ell},\what{H}_{\sg_2-\a_\ell}$ are linearly independent, as their images are so in $\hh$. 
			Also
			$\what{H}_{\sg_1+\sg_2-{\a_\ell}}=\what{H}_{\sg_1-\a_\ell}+\what{H}_{\sg_2-\a_\ell}+\what{H}_{\a_\ell}.$ Thus $\dim\what{\hh}=\dim\hh.$ The remaining types can be treated in a similar manner.\qed
				
			{  
				\begin{lem}\label{lemrev1}
					Let $\what Y=[\what{Y}_1,\ldots,\what{Y}_n]$, where $n\geq 2$ and $\what{Y}_i\in\{\what{X}_\a,\what{H}_\a,\what{d}_1,\what{d}_2\mid\a\in\pp^\pm\}$.
					Then $\what{Y}$ is in the span of brackets of the form $[\what X_{\a_1},\ldots,\what X_{\a_m}]$,
					$\a_i\in\pp^\pm$, $m\geq 1$.
				\end{lem}
			}
			
			{
				\proof
				We use induction on $n$. If $n=2$, then we are done using defining relations (\ref{eq111}).
				Assume now that $n>2$. By induction steps, we may assume without loss of generality that $\what Y=[\what{Y}_1,\what{Y}']$, where 
				$\what{Y}'=[\what X_{\a_1},\ldots,\what X_{\a_m}]$ for some $\a_i\in\pp^\pm$. Now if $\what{Y}_1=\hat{X}_\a$ for some $\a\in\pp^\pm$, then we are done. Otherwise, we have
				$$\what{Y}=[\what{Y}_1,[\what{X}_{\a_1},\what{X}']]=-[\what{X}',[\what{Y}_1,\what{X}_{\a_1}]]-[\what{X}_{\a_1},[\what{X}',\what{Y}_1]],$$
				with $\what{X}'=[\what {X}_{\a_2},\ldots,\what X_{\a_m}],$ and so we are done again by induction steps and by defining relations (\ref{eq111}).\qed
			}
			
			Using Lemma \ref{llem16}, we identify $\what\hh$ with $\hh$. For $\a\in\hh^\star$, we set
			$$\what\ll_\a:=\{\what X\in\what\ll\mid[h,\what X]=\a(h)\what X\hbox{ for all }h\in\hh\}.$$
			
			\begin{pro}\label{proper}
				The Lie algebra $(\what\ll,\what\hh)$ is an elliptic Lie algebra with root system $R$.
			\end{pro}
			
			\proof 
			{We first prove that $\what\ll=\sum_{\a\in\hh^\star}\what\ll_\a$. For this we need to show that any bracket in $\what\ll$ of the form $\hat Y=[\what Y_1,\ldots,\what Y_n]$, where $\what Y_i$'s belong to defining generators of $\what\ll$, sits in $\what\ll_\a$ for some $\a\in\hh^\star$.
				By defining relations of $\what\ll$, we have
				$\what X_\a\in\what\ll_\a$ and $\what H_\a, \what d_1,\what d_2\in\what\ll_0$, for $\a\in\pp^\pm$. So we are done if $n=1$. Assume now that $n\geq 2$. By Lemma \ref{lemrev1}, we may assume $\what Y=[\what X_{\a_1},\ldots,\what X_{\a_n}]$ for some $\a_i\in\pp^\pm$.
				Now if $\beta\in\pp^\pm$, then by the Jacobi identity and defining relations of $\what\ll$, we get
				$[\what H_\beta,[\what X_{\a_1},\what X_{\a_2}]]=(\a_1+\a_2)(\what H_\beta)[\what X_1,\what X_2].$
				Also  $[\what d_i,[\what X_{\a_1},\what X_{\a_2}]]=\lam_i(\a_1+\a_2) [\what X_1,\what X_2]=(\a_1+\a_2)(\what d_1)[\what X_{\a_1},\what X_{\a_2}].$ Thus
				$[\what X_{\a_1},\what X_{\a_2}]\in\what\ll_{\a_1+\a_2}.$
				Then an induction argument on $n\geq 2$ proves that
				$\what Y\in\what\ll_{\a_1+\cdots+\a_n}.$ The above reasoning shows that $\what\ll=\sum_{\a\in\hh^\star}\what\ll_\a$, and moreover, 
				$\what\ll_\a$ is spanned by brackets $[\what X_{\a_1},\ldots,\what X_{\a_n}]$ with
				$\a_1+\cdots+\a_n=\a$ where $\what X_{\a_i}\in\{\what X_\a,\what H_\a, \what d_1,\what d_t\mid\a\in\pp^\pm\}$.  This completes the proof that $\what\ll=\sum_{\a\in H^\star}\what\ll_\a$. Now defining relations of the form (V) implies that $\what\ll_\a=\sum_{\a\in R}\what\ll_\a$. Since $\Psi(\what\ll_\a)=\ll_\a$, we get $\what\ll_\a\not=\{0\}$ for $\a\in R$.\qed}

		
	\subsection{Type $A_1$}\label{A_1}
We specialize the presented Lie algebra $\what\ll$ to the case $X=A_1$ as an example. Since $\nu=2$, we have  up to isomorphism  two
elliptic root systems of type $A_1$:
$$R=\left\{\begin{array}{ll}
(\bbbz\sg_1\oplus\bbbz\sg_2)\cup(\pm\a+\bbbz\sg_1+2\bbbz\sg_2)
\cup(\pm\a+2\bbbz\sg_1+\bbbz\sg_2)&\hbox{if }\ind(R)=0,\\
(\bbbz\sg_1\oplus\bbbz\sg_2)\cup(\pm\a+\bbbz\sg_1+\bbbz\sg_2)&\hbox{if }\ind(R)=1.
\end{array}\right.$$ 
We fix the reflectable base
\begin{equation}\label{cc1}
\pp=\left\{\begin{array}{ll}
	\{\a,\sg_1-\a,\sg_2-\a\}&\hbox{if }\ind(R)=0,\\
	
	\{\a,\sg_1-\a,\sg_2-\a,\sg_1+\sg_2-\a\}&\hbox{if }\ind(R)=1,\end{array}
\right.
\end{equation}
for $R$.

			For a non-isotropic root $\be$, we define the automorphism $\Phi_\be\in\Aut(\what\ll)$ by
			$$\Phi_\be:=\exp(\ad\what{X}_\be)\exp(-\ad\what{X}_{-\be})
			\exp(\ad\what{X}_\be).$$

			\begin{lem}\label{lem234}
				$\Phi_{\sg_1-\a}(\what{X}_{\sg_2-\a})\in\what\ll_{\sg_2-2\sg_1+\a}$
				if and only if
				\begin{equation}
					\label{mic1}
					\frac{1}{2}\ad \what{X}_{\sg_1-\a}\ad \what{X}_{\a-\sg_1}[\what{X}_{\a-\sg_1},\what{X}_{\sg_2-\a}]=[\what{X}_{\a-\sg_1},\what{X}_{\sg_2-\a}].
				\end{equation}
			\end{lem}
			
			\proof We have
			\small{\begin{eqnarray*}
					&\Phi_{\sg_1-\a}(\what{X}_{\sg_2-\a})&\\
					&=&\hspace{-1cm}\exp\ad(\what{X}_{\sg_1-\a})(\what{X}_{\sg_2-\a}-[\what{X}_{\a-\sg_1},\what{X}_{\sg_2-\a}]{+}\frac{1}{2}[\what{X}_{\a-\sg_1},[\what{X}_{\a-\sg_1},\what{X}_{\sg_2-\a}]])\\
					&=&\hspace{-1cm}\what{X}_{\sg_2-\a}-[\what{X}_{\a-\sg_1},\what{X}_{\sg_2-\a}]-[\what{X}_{\sg_1-\a},[\what{X}_{\a-\sg_1},\what{X}_{\sg_2-\a}]]\\
					&&
					\hspace{-5mm}+\frac{1}{2}	[\what{X}_{\a-\sg_1},[\what{X}_{\a-\sg_1},\what{X}_{\sg_2-\a}]]
					+\frac{1}{2}
					[\what{X}_{\sg_1-\a},[\what{X}_{\a-\sg_1},[\what{X}_{\a-\sg_1},\what{X}_{\sg_2-\a}]]]\\
					&&\hspace{-.5mm}+\frac{1}{4}	[\what{X}_{\sg_1-\a},[\what{X}_{\sg_1-\a},[\what{X}_{\a-\sg_1},[\what{X}_{\a-\sg_1},\what{X}_{\sg_2-\a}]]]].
				\end{eqnarray*}
			}
			We note that we have
			\small{\begin{eqnarray*}
					[\what{X}_{\sg_1-\a},[\what{X}_{\a-\sg_1},\what{X}_{\sg_2-\a}]]&=&
					-[\what{X}_{\sg_2-\a}[\what{X}_{\sg_1-\a},\what{X}_{\a-\sg_1}]]-
					[\what{X}_{\a-\sg_1},[\what{X}_{\sg_2-\a},\what{X}_{\sg_1-\a}]]\\
					&=&
					-[\what{X}_{\sg_2-\a},H_{\sg_1-\a}]=(\a,\a^\vee)\what{X}_{\sg_2-\a}=2\what{X}_{\sg_2-\a}.
				\end{eqnarray*}
			}
			Then from this and  (\ref{mic1}), we get
			$$\begin{array}{l}
				\what{X}_{\sg_2-\a}-[\what{X}_{\sg_1-\a},[\what{X}_{\a-\sg_1},\what{X}_{\sg_2-\a}]]+\frac{1}{4}	[\what{X}_{\sg_1-\a},[\what{X}_{\sg_1-\a},[\what{X}_{\a-\sg_1},[\what{X}_{\a-\sg_1},\what{X}_{\sg_2-\a}]]]]\vspace{2mm}\\
				\hspace{1cm}=\what{X}_{\sg_2-\a}-2\what{X}_{\sg_2-\a}+\what{X}_{\sg_2-\a}=0,
			\end{array}
			$$
			and
			\small{
				\begin{eqnarray*}
					-[\what{X}_{\a-\sg_1},\what{X}_{\sg_2-\a}]
					&+&\frac{1}{2}
					[\what{X}_{\sg_1-\a},[\what{X}_{\a-\sg_1},[\what{X}_{\a-\sg_1},\what{X}_{\sg_2-\a}]]]\\
					&=&
					-[\what{X}_{\a-\sg_1},\what{X}_{\sg_2-\a}]+[\what{X}_{\a-\sg_1},\what{X}_{\sg_2-\a}]=0.
				\end{eqnarray*}
			}
			\qed
			
			\begin{lem}\label{canc1}
				Suppose $R$ is an elliptic root system of type $A_1$ and $\pp$ is the reflectable base given in (\ref{cc1}) for $R$. Let $\what\ll$ be the Lie algebra defined by generators  
				$\what{X}_{\a},\;\what{H}_\a,\what d_1,\what d_2$, $\a\in\pp^\pm$, and relations I-VI given in (\ref{eq111}) plus relation (\ref{mic1}).
				Then $(\what\ll,\what\hh)$ is an elliptic Lie algebra satisfying
				$$
				\Phi_{\gamma}(\what\ll_\be)=\what\ll_{w_\gamma(\be)},\qquad(\gamma\in\pp,\be\in R^\times).
				$$ 
			\end{lem}
			
			\proof We first note that (\ref{mic1}) holds in $\ll$ for $e_{\sg_1-\a}, e_{\sg_2-\a}$ in place of $\what X_{\sg_1-\a}$ and $\what  X_{\sg_2-\a}$. Thus the epimorphism $\Psi:\what\ll\rightarrow\ll$ given in (\ref{eq1234})  holds in our case. Suppose $\gamma,\be$ are as in the statement.
			Let $\what X\in\what\ll_\be$. We may assume without loss of generality 
			that $\what X=[\what{X}_{\be_1},[\what{X}_{\be_2},[\ldots,[\what{X}_{\be_{k-1}},\what{X}_{\be_k}]\ldots]$ for some $\be_i\in\pp$ with $\be=\be_1+\cdots+\be_k$.
			Then for $\gamma\in\pp$, we have
			$$
			\Phi_\gamma(\what X)=
			[\Phi_{\gamma}(\what{X}_{\be_1}),[\Phi_{\gamma}(\what{X}_{\be_2}),[\ldots,[\Phi_{\gamma}(\what{X}_{\be_{k-1}}),\Phi_{\gamma}(\what{X}_{\be_k})]\cdots]].
			$$
			Now if $\{\gamma,\be_i\}\not=\{\sg_1-\a,\sg_2-\a\}$, then the vectors
			$\what{X}_{\gamma},\what{X}_{\be_i}$ can be identified as root vectors of a subalgebra of $\what\ll$ which is either finite-dimensional simple, or {derived subalgebra of an affine Lie algebra}, implying that
			$\Phi_\gamma(\what{X}_{\be_i})\in\what\ll_{w_{\gamma(\be_i)}}$.
			If $\{\gamma,\be_i\}=\{\sg_1-\a,\sg_2-\a\}$, then
			since (\ref{mic1}) holds by assumption, we get from Lemma \ref{lem234} that $\Phi_\gamma(\what{X}_{\be_i})\in\what\ll_{w_{\gamma(\be_i)}}$.
			Thus $\Phi_\gamma(\what\ll_\be)\sub\what\ll_{w_\gamma(\be)}$.
			This completes the proof.\qed
			
			\subsection{Type $A_1$, $\ind(R)=0$}\label{pg1}
			We now restrict our attention to the case when $R$ is of type $A_1$ and $\ind(R)=0$. We consider the corresponding reflectable base  
			$\pp$ given in (\ref{cc1}). Let $\what\ll$ be the corresponding presented Lie algebra. 
		
	We recall that for any elements
${X}_0,\ldots {X}_m$ of a Lie algebra and for $1\leq j<m$, we have
\begin{equation}\label{ten}
\begin{array}{lcl} 	[{X}_0,\ldots {X}_m]&=&[[{X}_0,{X}_1],{X}_2,\ldots,{X}_m]\\ &+&[{X}_1,[{X}_0,{X}_2],{X}_3,\ldots,{X}_m]\\
	&+&\cdots\\ 	&+&[{X}_1,{X}_2,\ldots, {X}_{j-1}, [{X}_0,{X}_{j}],{X}_{j+1},\ldots {X}_m]\\
	&+&[{X}_1,\ldots, {X}_{j-1}, {X}_0, {X}_j, \ldots, {X}_m].
\end{array}
\end{equation}
		
		\begin{lem}\label{lem235}
			{Let $\be,\a_1,\ldots,\a_m\in\pp^\pm$ with $\a_1+\cdots+\a_m=\be$. Then}
			
			(i) $m$ is odd, and for each $\gamma\in\pp^{\pm}$ with $\gamma\not=\be$, the number of $j$'s such that $\a_j=\gamma$, is even. 
			
			(ii) If
			$A=[\what{X}_{\a_1}, \what{X}_{\a_2},\ldots,\what{X}_{\a_m}]\in{\what\ll}_\be$, then $A$ can be written as a sum of expressions of the form
			$[\what{X}_{\be_1},\ldots,\what{X}_{\be_k}]$, where all $\be_i$'s belong to
			$\pm\{\a,\sg_1-\a\}$ or all belong to $\pm\{\a,\sg_2-\a\}$.  	
		\end{lem}
		
		\proof (i) is clear as $\a,\sg_1,\sg_2$ are linearly independent. 
		
		(ii)  We recall that 
		$R=R^0\cup R^\times$ with
		$$R^0=\bbbz\sg_1\oplus\bbbz\sg_2\andd R^\times
		=(\pm\a+2\bbbz\sg_1\oplus\bbbz\sg_2)\cup
		(\pm\a+\bbbz\sg_1\oplus2\bbbz\sg_2).$$
		In particular, we have
		\begin{equation}\label{temp234}
			\ep\a+n_1\sg_1+n_2\sg_2\not\in R\hbox{ for } \ep\in\{0,\pm1\},\; n_1,n_2\in2\bbbz+1.
		\end{equation}
		This fact will be used frequently in the proof without further mention. We also recall that $\pp=\{\a,\sg_1-\a,\sg_2-\a\}$.
		
		Without loss of generality, we may assume $\be=\a$. If all $\a_i$'s belong to $\{\a,\sg_i-\a\}^\pm$, $i=1,2$, we are done. So we may assume that
		\begin{equation}\label{temp235}
			\a_i=\pm(\sg_1-\a)\hbox{ and }
			\a_j=\pm(\sg_2-\a)\hbox{ for some }i,j.
		\end{equation}
		If $m=1$, there is nothing to prove. Now (\ref{temp234}), (\ref{temp235}) {and relations of the form (V$'$) show that no nonzero expression} of the form $A$ exists for length $m=3$. By part (i), for $m=5$ the only possibility is $\{\a_1,\ldots,\a_5\}=\{\a,\pm(\sg_1-\a),\pm(\sg_2-\a)\}
		$ and the only possibility for $A$, up to a permutation of $\pp$, is $[\what{X}_{\a-\sg_2},\what{X}_{\sg_2-\a},\what{X}_{\a-\sg_1},\what{X}_{\sg_1-\a},\what{X}_\a]$ which reduces to the case $m=3$, by using (\ref{ten}) and applying defining relations of the form (III). In general
		by using (\ref{ten}), we write $A$ as a sum of terms of length $m-2$, except for the last term
		$[\what{X}_1,\ldots,\what{X}_{j-1},\what{X}_0,\what{X}_j,\ldots,\what{X}_m]$ which has length $m$. But then using (\ref{temp234}), we choose an appropriate $j$ such that this term becomes zero, {by applying a relation of the form (V$'$)}.
		\qed

		\begin{thm}\label{tmba1}
			Let $R=\Lam\cup(\pm\a+S)$ be an elliptic root system of type $A_1$ where $\Lam=\bbbz\sg_1\oplus\bbbz\sg_2$ and
			$S=(\sg_1+2\Lam)\cup(\sg_2+2\Lam)$. Let $(\what\ll,\what\hh)$ be the Lie algebra defined by generators and relations (\ref{eq111}).
			Then $\what\ll$ is an elliptic Lie algebra with root system $R$. Moreover $\dim\what\hh=5$ and 
			$\dim{\what\ll}_\be=1$ for $\be\in R^\times$.
		\end{thm} 
		
	\proof By Proposition \ref{proper} and Lemma \ref{llem16}, $\what\ll$ is an elliptic Lie algebra with root system $R$ satisfying $\dim\what\hh=5$. Recall from Proposition \ref{coron2} that if
	$\pp=\{\a,\sg_1-\a\}$ or $\pp=\{\a,\sg_2-\a\}$, then the root system $R_\pp\sub R$ is an affine root system of type $A_1$. 
	We now consider the Lie algebra $\what\ll_\pp$ defined by generators
	$\what X_{\a}$, $\what H_\a$, $\a\in\pp^{\pm}$, subject to relations (\ref{eq111})(I)-(V). {By Remark \ref{revrem2}}, we see that $\what\ll_\pp$ is isomorphic to the derived algebra of an affine Lie algebra of type $A_1$. Thus identifying $\what\ll_\pp$ as a subalgebra of $\what\ll$, and using Lemma \ref{lem235}, we conclude that $\dim{\what\ll}_\be=1$ for $\be\in\pp$. We note that by (\ref{temp234}) and (\ref{eq111}(V)), both sides of equation (\ref{mic1}) are zero and so this relation is a consequence of relations holds in $\what\ll$.
	Thus Lemma \ref{canc1} implies that 
	$\dim{\what\ll}_\be=1$ for all $\be\in R^\times$. 
	\qed
	
	\subsection{Further considerations}\label{fur}
	In \cite{SaY, Ya}, the authors construct elliptic Lie algebras using Serre's type generators and relations. {The defining generators of the given presentations are based on a} ``root base'' assigned to the corresponding elliptic root system. This root base which we call a {\it Saito root base} and denote it by $\pp_S$ is defined by K. Saito \cite[Section (5.2)]{Sa}.
	For a given type, in Table \ref{tabnew} the cardinality of a Saito root base is compared with the corresponding reflectable base of Table \ref{ttt1}.

		{\scriptsize
			\begin{table}[ht]
				\caption{Cardinals of Saito bases and reflectable bases}
				\label{tabnew}
				\vspace{-4mm}
				\begin{tabular}
					[c]{|c |   c |c|} \hline
					$\mbox{Type } X$& $|\pp_S(X)|$ &$|\pp(X)|$ \\
					\whline
					$A _1$ & $\begin{array}{c}
						$3$\hbox{ if }\ind(R)=0\\
						$4$\hbox{ if }\ind(R)=1
					\end{array}$&$\begin{array}{c}
						$3$\hbox{ if }\ind(R)=0\\
						$4$\hbox{ if }\ind(R)=1
					\end{array}$\\
					\hline
					$A_\ell (\ell>1)$&$2\ell+2$&$\ell+2$\\
					\hline
					$D_\ell$&$2\ell-2$&$\ell+2$\\
					\hline
					$E_{6,7,8}$ & $8,9,10$&$8,9,10$\\
					\hline
					$F_4$ & $\begin{array}{c}
						$6$\hbox{ if }t=0\\
						$6$\hbox{ if }t=1\\
						6\hbox{ if }t=2
					\end{array}
					$& $\begin{array}{c}
						$6$\hbox{ if }t=0\\
						$6$\hbox{ if }t=1\\
						6\hbox{ if }t=2
					\end{array}
					$\\
					\hline
					$G_2$ & $\begin{array}{c}
						4\hbox{ if }t=0\\
						4\hbox{ if }t=1\\
						4\hbox{ if }t=2
					\end{array}
					$& $\begin{array}{c}
						4\hbox{ if }t=0\\
						4\hbox{ if }t=1\\
						4\hbox{ if }t=2
					\end{array}
					$\\
					\hline
					$B_2$ & $\begin{array}{c}$5$\hbox{ if }t=0,\;\ind(R)=0\\
						6\hbox{ if }t=0,\;\ind(R)=1\\	
						4\hbox{ if }t=1
						\\
						5\hbox{ if }t=2,\;\ind(R)=0\\
						6\hbox{ if }t=2,\;\ind(R)=1
					\end{array}
					$&$\begin{array}{c}4\hbox{ if }t=0,\;\ind(R)=0\\
						5\hbox{ if }t=0,\;\ind(R)=1\\	
						4\hbox{ if }t=1
						\\
						{4}\hbox{ if }t=2,\;\ind(R)=0\\
						5\hbox{ if }t=2,\;\ind(R)=1
					\end{array}
					$\\
					\hline
					$B_\ell\; (\ell>2)$ & $\begin{array}{c}
						2\ell-1\hbox{ if }t=0\\
						2\ell\hbox{ if }t=1
						\\
						2\ell+1\hbox{ if }t=2,\;\ind(R)=0\\
						2\ell+2	\hbox{ if }t=2,\;\ind(R)=1
					\end{array}
					$& $\begin{array}{c}
						\ell+2\hbox{ if }t=0\\
						\ell+2\hbox{ if }t=1
						\\
						\ell+2\hbox{ if }t=2,\;\ind(R)=0\\
						\ell+3	\hbox{ if }t=2,\;\ind(R)=1
					\end{array}
					$\\
					\hline
					$C_\ell\; (\ell>2)$ & $\begin{array}{c}2\ell+2\hbox{ if }t=2\\
						2\ell\hbox{ if }t=1
						\\
						2\ell+1\hbox{ if }t=0,\;\ind(R)=0\\
						2\ell+2\hbox{ if }t=0,\;\ind(R)=1
					\end{array}
					$&$\begin{array}{c}\ell+2\hbox{ if }t=2\\
						\ell+2\hbox{ if }t=1
						\\
						\ell+2\hbox{ if }t=0,\;\ind(R)=0\\
						\ell+3\hbox{ if }t=0,\;\ind(R)=1
					\end{array}
					$\\
					\hline
				\end{tabular}
			\end{table}   
		}
		
		{Let $R$ be an elliptic root system of rank $>1$. Let $\pp$ be the reflectable base for $R$ given in Table \ref{ttt1}, and $\pp_S$ be the Saito base for $R$. We have $|\pp|\leq |\pp_S|$, see Table \ref{tabnew}. Considering some justifications, we may identify the reflectable base $\pp$ with a subset of $\pp_S$. Let $\what\ll_S$ be the corresponding presented Lie algebra given in \cite{SaY, Ya} associated to $\pp_S$. 	 Consider the map which assigns each generator of $\what\ll$ to the corresponding generator in $\what\ll_S$, see \cite[Definition 2]{SaY} and \cite[Definition 3.1]{Ya}. One sees from \cite[Theorem 2]{SaY} and \cite[Definition 3.1 and Theorem 3.1]{Ya} that the defining relations (I)-(IV) and (VI) of (\ref{eq111}) hold in $\what\ll_S$ for the corresponding generators. 
			Moreover, since $R$ is the set of roots of $\what\ll_S$, it follows that
			the relations of the form (V) also hold in $\what\ll_S$. Therefore, we get a
			homomorphism $\Phi:\what\ll\rightarrow\what\ll_S$ which preserves root spaces. This suggests that to achieve a finite Serre's type presentation associated with $\pp$, with $1$-dimensional (non-isotropic) root spaces, one should investigate replacing relations of the form (V) with more appropriate relations such as (V$'$) and those given in \cite{SaY} and \cite{Ya}.} 

	\end{document}